\documentclass[twoside]{article}
\newlength{\ei}\ei=0.0138888889em   
   \newlength{\breit}
\newlength{\SyW}  \newlength{\msu}  \msu=\mathsurround 
\newcommand{\weg}{\hspace{\msu}}\newcommand{\ran}{\hspace{-\mathsurround}}
\newcommand{\qed}{\nopagebreak\hspace*{\fill}\mbox{\scriptsize\rm
  {/}\hspace{-13\ei}{/}\hspace{-13\ei}{/}\hspace{-13\ei}{/}}}   
\newcommand{\nz}{\normalsize}\newcommand{\sz}{\small}

\newcommand{\Ds}{\displaystyle} \newcommand{\Ts}{\textstyle}
  
\newcommand{\so}{\;\Longrightarrow\;}

\newcommand{\Nsup}{\sup\nolimits}
\newcommand{\Ninf}{\inf\nolimits}

\newcommand{\Nlim}{\lim\nolimits}
\newcommand{\Nlimsup}{\limsup\nolimits}
\newcommand{\Nliminf}{\liminf\nolimits}
\newcommand{\Tfrac}[2]{{\Ts\frac{#1}{#2}}}

\newcommand{\ave}{\mathop{\rm ave\hskip6\ei}\nolimits}
\newcommand{\rave}{\mathop{\Ts\sqrt{n}\hskip12\ei\rm ave\hskip6\ei}\nolimits}
\newcommand{\itg}{\int}

\newcommand{\ggrto}{\relbar\joinrel\longrightarrow}

\newcommand{\gwto}{\mathrel{\mathsurround0em \mbox{$\ggrto$}%
  \llap{\settowidth{\SyW}{$\ggrto$}
  \raisebox{-.15ex}{\makebox[\SyW]{\scriptsize\rm w}}}}}
\newfont{\tenmsam}{msam10}                                 

\font\hMtenbf=cmbx10 scaled \magstephalf 

\newcommand{\zi}{\mathop{\rm sign}\nolimits}
\newcommand{\Jc}{\mathop{\bf {{}I{}}}\nolimits}
\newcommand{\hJc}{\mathop{\mbox{\hMtenbf I}}\nolimits}

\newcommand{\Lo}{\mathop{\rm {{}o{}}}\nolimits} 
\newcommand{\bvLo}{\mathop{\rm {{}\breve{o}{}}}\nolimits}
\newcommand{\tlLo}{\mathop{\rm {{}\tilde{o}{}}}\nolimits}
\newcommand{\htLo}{\mathop{\rm {{}\hat{o}{}}}\nolimits}
\newcommand{\brLo}{\mathop{\rm {{}\bar{o}{}}}\nolimits}
\newcommand{\lin}{\mathop{\rm lin}\nolimits}
\newcommand{\cl}{\mathop{\it c\ell}\nolimits}
\newcommand{\clin}{\mathop{\cl\lin}\nolimits}

\newcommand{\Ew}{\mathop{\rm {{}E{}}}\nolimits} 

\makeatletter                  
\def\magstepams#1{\ifcase#1 \magstep0\or\magstephalf\or\magstep1\fi\relax}
\font\tenmsbm=msbm10 scaled \magstepams\@ptsize
\font\sevmsbm=msbm7  scaled \magstepams\@ptsize
\def\relaxnext@{\let\next\relax}   \def\noaccents@{\def\accentfam@{0}}
\newfam\msbmfam \textfont\msbmfam=\tenmsbm \scriptfont\msbmfam=\sevmsbm
\def\Bberr{Use Black board bold only in math mode}
\def\Bb{\relaxnext@\ifmmode\let\next\Bb@\else
        \def\next{\errmessage{\Bberr}}\fi\next}
\def\Bb@#1{{\Bb@@{#1}}}\def\Bb@@#1{\noaccents@\fam\msbmfam#1}\makeatother
\newcommand{\R}{{\mbox{\tenmsbm R}}}        
\newtheorem{Thm}{Theorem}[section]   
\newtheorem{Prop}[Thm]{Proposition}  \newtheorem{Lem}[Thm]{Lemma}
     
\newtheorem{Exa}[Thm]{Example}       \newtheorem{Rem}[Thm]{Remark}
\newenvironment{Bew}{\begin{trivlist}\item[]{\sc Proof}}{\end{trivlist}}

\newcommand{\rfi}[1]{\makebox[\parindent][l]{%
                     \makebox[0em][r]{\rm(}\sf#1\rm)}}
\newcounter{ABCc}\renewcommand{\theABCc}{\alph{ABCc}}
\newenvironment{ABC}{\begin{list}{
  \rfi{\theABCc}}{\usecounter{ABCc} \topsep 0ex \partopsep 0ex \itemsep0ex
  \parsep=\parskip \leftmargin 0em \rightmargin 0em \itemindent=\parindent
  \listparindent=\parindent  \labelsep 0em \labelwidth 0em }}{\end{list}}
\newcounter{ABCca}\renewcommand{\theABCca}{\arabic{ABCca}}

\begin{document}
\title{\LARGE\bf One-Sided Confidence About Functionals Over Tangent Cones}
\author{\nz\rm Helmut Rieder\\
        \sz\sl University of Bayreuth, Germany}
\date{\sz\rm 1.~March~2002, under revision}
\maketitle
\begin{abstract}\noindent
In the setup of i.i.d.~observations and a real valued differentiable
functional~$T$, locally asymptotic upper bounds are derived for the
power of one-sided tests (simple, versus large values of~$T$)
and for the confidence probability of lower confidence limits
(for the value of~$T$), in the case that the tangent set is only
a convex cone.
The bounds, and the tests and estimators which achieve the bounds,
are based on the projection of the influence curve of the functional
on the closed convex cone, as opposed to its closed linear span.
The higher efficiency comes along with some weaker, only one-sided,
regularity and stability.
\par \noindent
{\sl AMS/MSC-2000 classification}: 62F35.
\end{abstract}
\section{Introduction}                \setcounter{equation}{0}\label{s.I}
Given a model~${\cal P}$ of probability measures on some sample space,
let some one dimensional aspect be defined by some statistical functional
$T \colon {\cal P}\to\R$. We consider the simplest case of~$n$
stochastically independent observations $ x_1,\ldots,x_n $ with identical
distribution any $P\in {\cal P}$, and the task is to make confidence
statements on the unknown value~$T(P)$ by means of tests and estimators.
\par                                        
In the usual testing problems concerning the value of~$T$, the power of
level~$\alpha$ tests cannot exceed certain asymptotic upper bounds.
Likewise, the accuracy of estimators of~$T(P)$ is limited by some asymptotic
  upper bounds for one- and two-sided confidence probabilities.
These bounds form a classical subject of non- and semiparametric theory;
confer, for example, Bickel et al.~(1993), Pfanzagl and Wefelmeyer~(1982),
Rieder~(1994), and van der Vaart~(1998).
\par
Having fixed any $P\in {\cal P}$, either for the purpose of testing local
alternatives or, in estimation, to be able to exclude artificial phenomena
of superefficiency, local variations of~$P$ within~${\cal P}$ must be taken
into account\footnote{implicitly, already, in the classical scores
function---a derivative, of log densities.}.
These variations are formulated as differentiable
paths~${(P_{g,s})}_{s>0}$ in~${\cal P}$, in direction of
certain tangents~$g\in \nolinebreak L_2(P)$ at~$P$,
such that, in the Hilbert space of square root densities,
\begin{equation} \label{e:i:pathdiff}
  \sqrt{dP_{g,s}}\,=
  \bigl(1+ \Tfrac{1}{2}s \hspace{6\ei} g \bigr) \sqrt{dP}\,
  + \Lo(s) \qquad   \mbox{as $\Ds s \downarrow 0$} \hspace{-1em}
\end{equation}
The functions~$g$ necessarily have expectation
  $\Ew g=\langle g|1\rangle=0$ under~$P$; in other words,
  $g\perp \mbox{the constants}$ in $L_2(P)$.
Given any $g\in L_2(P)$, $\langle g|1\rangle=0$,
a corresponding path (in the set of all probabilities) is
\begingroup \mathsurround0em\arraycolsep0em \begin{eqnarray}
\label{E:I:pathdef} dP_{g,s} & {}={} &
  \Bigl(\Tfrac{1}{2}s \hspace{6\ei} g + \sqrt{
  1- \Tfrac{1}{4} s^2 \hspace{6\ei}\Vert g\Vert^2 \hspace{6\ei}}
  \hspace{24\ei}\Bigr)^2 \,dP \\
\noalign{\noindent or\nopagebreak} \label{e.i.dPs1fach}
   dP_{g,s} & {}={} & (1+s \hspace{6\ei} g)\,dP  \hspace{2em}
   \mbox{\mathsurround\msu if $\Ds g\in L_{\infty}(P) $}
\end{eqnarray}\endgroup
The set~${\cal G}$ of all tangents at~$P$ on one hand reflects the richness
of the model~${\cal P}$. On the other hand, ${\cal G}$~is restricted by the
differentiability requirement on the functional:
There exist some function $\kappa\in L_2(P)$, such that for every
$g\in {\cal G}$ and any path~(\ref{e:i:pathdiff}) in~${\cal P}$,
\begin{equation} \label{e:i:Tdiff}
  T(P_{g,s})=T(P)+s\langle \kappa|g\rangle + \Lo(s) \qquad
  \mbox{as $\Ds s \downarrow 0$} \hspace{-.75em}
\end{equation}
The function $\kappa$, a so-called influence curve of~$T$ at~$P$, may
not be unique. But the orthogonal projection~$\bar{\kappa}$ of~$\kappa$
on the closed linear span~$\clin{{\cal G}}$ of~${\cal G}$ in~$L_2(P)$
is unique---the canonical gradient, or efficient influence curve.
\par     By definition, the tangent set~${\cal G}$ of~${\cal P}$ at~$P$
is a cone in~$L_2(P)\cap \{1\}^{\perp}$ with vertex at~$0$,
such that $\gamma g\in{\cal G}$ for $g\in {\cal G}$ and
          $\gamma\in [\hskip6\ei0,\infty)$.
For example, the classical nonparametric alternative hypotheses of
positive asymmetry and positive dependence naturally lead to cones.
\par
Furthermore, there is a general argument why arbitrary tangent sets
should be considered in theory. In testing, the null hypothesis usually
is canonical and simple, but the alternative may be chosen freely,
more complex, according to the particular case at hand.
In estimation, as noted by one referee, tangent cones arise if the
paramater value is a boundary point of the parameter set.
Moreover, also for other parameter values, the previous argument
may be enforced from a robustness viewpoint.
In the setup of Rieder~(1994; Chapter~4), any parametric model
distribution may be enlarged to infinitesimal neighborhoods consisting
of the local alternatives generated by, for example, a tangent cone
(leading us to consider the smallest cone containing the neighborhood
cone and the linear span of the parametric tangent).
\par
In most of the literature on asymptotic bounds so far, the tangent
set is assumed a linear space~${\cal G}=\bar{{\cal G}}$, such that
  $\clin \bar{{\cal G}}$ is just the closure~$\cl \bar{{\cal G}}$
of~$\bar{{\cal G}}$. Then the said bounds are determined by the canonical
gradient~$\bar{\kappa}$, acting as a least favorable (limiting) tangent,
and its norm~$\Vert \bar{\kappa}\Vert $. 
\par
If the tangent set is not a linear space but a possibly nonconvex
cone~${\cal G}=\tilde{{\cal G}}$, the situation is not quite
clear\footnote{As for nonconvex cones, we refer to the footnote summary
in van der Vaart~(1998; p~367).}.
In our paper, we shall settle on cones~$\tilde{{\cal G}}$ that in
addition are convex, such that
  $ \gamma_1 g_1+ \gamma_2 \hspace{3\ei}g_2 \in \tilde{{\cal G}} $
for $ g_i\in \tilde{\cal G}$ and $ \gamma_i\in [\,0,\infty) $.
\par
Even in this case, of a convex tangent cone~$\tilde{{\cal G}}$,
the results in literature seem somewhat contradictory:
On~one hand, the convolution representation and asymptotic minimax risk
under symmetric subconvex loss given by van der Vaart~(1998; Theorems~25.20
and~25.21) are still expressed by the canonical gradient~$\bar{\kappa}$
(the orthogonal projection of~$\kappa$ on~$\clin \tilde{{\cal G}}$).
On~the other hand, Pfanzagl and Wefelmeyer~(1982; Theorem~9.2.2)
state a two-sided concentration bound in terms of the (smaller)
projection~$\tilde{\kappa}$ of~$\kappa$ on a closed convex tangent
cone~$\tilde{{\cal G}} = \nolinebreak \cl \tilde{{\cal G}}$.
Their proof, however, makes use of
  $ - \tilde{{\cal G}} \subset \tilde{{\cal G}} $,
so their cone must in fact be a (closed) linear space.
Also Janssen~(1999), in the context of testing, considers convex
tangent cones~$\tilde{{\cal G}}$ and argues by the
projection~$\tilde{\kappa}$ of~$\kappa$ on~$\cl \tilde{{\cal G}}$.
But, throughout his paper, he treats~$\tilde{\kappa}$ as if it
were~$\bar{\kappa}$, as he nowhere accounts for the nonorthogonality
of the residual $\kappa-\tilde{\kappa}$ on~$\tilde{{\cal G}}$ 
in the case that $\bar{\kappa}\notin \cl\tilde{{\cal G}}$.
\par
Thus, either by result or by implicit assumption, the asymptotic power
and concentration bounds obtained so far for convex tangent cones agree
with those for their linear spans.
\par     The present investigation, in the case of convex tangent
cones~$\tilde{{\cal G}}$, derives locally asymptotic upper bounds
for the power of one-sided tests (of a simple hypothesis against large
values of~$T$), as well as for the confidence probabilities of lower
confidence limits for~$T(P)$. These asymptotic bounds
are given truly in terms of the projection~$\tilde{\kappa}$
of the influence curve~$\kappa$ of the functional~$T$
on the closed convex cone~$\cl\tilde{{\cal G}}$
(Theorems~\ref{t:t:1s-power} and~\ref{t:e:conf.bd}).
Since
  $\bar{\kappa}\in \clin \tilde{{\cal G}} \setminus \cl \tilde{{\cal G}}$
in general, that is, $ \tilde{\kappa} \ne \bar{\kappa} $ or,
equivalently, $\Vert \tilde{\kappa}\Vert < \Vert \bar{\kappa}\Vert $,
the upper bounds are larger than those based on~$\bar{\kappa}$.
\par
For the higher efficiency, however, a considerable price has to be paid,
which constists in a weaker and merely one-sided regularity and stability:
In the case of testing, the asymptotic size rises to~$100\%$ over an
only slightly enlarged, and therefore over the larger one-sided,
null hypothesis (Proposition~\ref{p.t.sizextH1}).
In~the case of estimation, the asymptotic bias may become plus
infinity under local alternatives (Proposition~\ref{p.cr.posmed}).
As a consequence, and as the (positive parts of) efficient estimators
are asymptotically unique
(Proposition~\ref{p.e.u}, Remark~\ref{r.cr.pfanz}),
the bound stated by Pfanzagl and Wefelmeyer~(1982; Theorem~9.2.2)
cannot possibly be attained under the condition of
asymptotic median unbiasedness.
The merely one-sided regularity, and one-sided asymmetric testing
pseudo-loss function, are also responsible for the 
difference to van der Vaart's~(1998) results.
\par
The investigation originated from the attempt by Rieder~(2000) to
\mbox{subject} robust statistics to the semiparametric approach by treating
neighborhoods as nuisance parameters, which leads to (the subtraction of the)
nonlinear projection on balls (from the classical scores).
Except for one-sided robust testing, however, the influence curves
thus obtained may differ from the optimally robust ones
of Rieder~(1994; Chapter~5).
Thus, contrary to what one would hopefully expect, the projection recipe
does not always give the optimal procedures. 
\par
Therefore, the present extension from linear spaces to convex cones
requires subtle modifictions of the proofs in the classical case.
Once derived, the new results ask for a careful interpretation of the
assumed regularity (Subsection~\ref{ss.e.reg}) and the implied stability
(Subsection~\ref{ss.e.lbh}), and a comparison for convex cones and their
linear spans becomes due (Subsection~\ref{ss.e.vs}).
\par
For reasons of comparability, throughout this paper the cases of
a linear tangent space~$\bar{{\cal G}}$ and a convex tangent
cone~$\tilde{{\cal G}}$, respectively, are stated together.
The $L_2(P)$-closure $\cl \bar{{\cal G}}$ of a linear tangent
space~$\bar{{\cal G}}$ is again a linear space, the $L_2(P)$-closure
    $\cl \tilde{{\cal G}}$ of a convex cone~$\tilde{{\cal G}}$ again
a convex cone. The canonical gradient, which is the projection of~$\kappa$
on~$\cl\bar{{\cal G}}$ and~$\clin \tilde{{\cal G}}$, respectively,
is denoted by~$\bar{\kappa}$, the projection of~$\kappa$
on~$\cl\tilde{{\cal G}}$ is denoted by~$\tilde{\kappa}$.
\par     
Convenient characterizations of the projections are supplied in the
appendix; the criteria (\ref{e:X:lin}) and~(\ref{e:X:cone}) for
$\bar{\kappa}$ and~$\tilde{\kappa}$ will be used without explicit reference.
Throughout the paper, the influence curve~$\kappa$, the tangent
space~$\bar{{\cal G}}$ and convex tangent cone~$\tilde{{\cal G}}$ at~$P$
are assumed of such a kind that
\begin{equation} \label{e.i.kkne0}
  \bar{\kappa}\ne0\weg, \qquad \tilde{\kappa}\ne0
\end{equation}
As noted, the interesting case occurs if 
  $ \bar{\kappa}\ne \tilde{\kappa} $.
\par     One-sided inference about non-smooth functionals of a density
has been studied by Donoho~(1988), by entirely different techniques and
in an even more nonparametric setting. 
Nevertheless, we encounter a somehow similar impossibility of
sensible upper confidence limits:
The estimators that provide the best lower confidence limits, subject to
some local asymptotic median nonnegativity, necessarily achieve overshoot
probability~$100\%$ 
under local alternatives. This distinguishes convex tangent cones from
linear tangent spaces, where the efficient estimator is unique and
asymptotically median unbiased.
\paragraph{Notation}
\ran $\langle \ldotp|\ldotp \rangle$
has already been used to denote the inner product in~$L_2(P)$. 
$\Jc$~stands for the indicator function.
Limits $\liminf_n$, $\limsup_n$, and~$\lim_n$ are meant for $n\to \infty$.
Asy.\ is our abbreviation of asymptotic/asymptotically.
\section{One-Sided Tests}             \setcounter{equation}{0}\label{s.T}
\subsection{Definition of Hypotheses}
For the fixed probability~$P\in {\cal P}$ and tangent set~${\cal G}$,
simple and one-sided composite asy.~hypotheses about the sequence of
laws~$Q_n$ of the i.i.d.~observations at sample size~$n=1,2,\ldots$
are defined by
\begingroup \mathsurround0em\arraycolsep0em
\begin{eqnarray}
\hspace{-1.5em}\label{e:t:J0}   J^0 &{}:{}& \hspace{30\ei}
                 Q_n=P \hspace{60\ei} \mbox{eventually}  \\
\hspace{-1.5em}\label{e:t:H0}   J &{}:{}& \hspace{30\ei} \Nlim_{n}
                 \sqrt{n}\,\bigl(\hspace{6\ei}T(Q_n)-T(P)\bigr) =0  \\
\hspace{-1.5em}\label{e:t:H1}   H &{}:{}& \hspace{30\ei}\Nlimsup_{n}
                  \sqrt{n}\,\bigl(\hspace{6\ei}T(Q_n)-T(P)\bigr)\le0  \\
\hspace{-1.5em}\label{e:t:K}    K &{}:{}& \hspace{30\ei}\Nliminf_{n}
                 \sqrt{n}\,\bigl(\hspace{6\ei}T(Q_n)-T(P)\bigr)\ge c
\end{eqnarray}\endgroup
where $c\in (0,\infty)$ is some fixed constant.
The measures~$Q_n$ in (\ref{e:t:H0})--(\ref{e:t:K}) may not be
arbitrary elements of model~${\cal P}$ but are assumed to approach~$P$
along any path~${(P_{g,s})}_{s>0}$ in~${\cal P}$ such that,
for some $g\in {\cal G}$ and $t\in (0,\infty)$, eventually,
\begin{equation} \label{e.t.Qpath}
    Q_n=P_{n,t,g}= P_{g,t/\!\sqrt{n}\,}
\end{equation}
In particular, every such sequence~$(Q_n^n)$ is contiguous to~$(P^n)$.
Also, the expansion~(\ref{e:i:Tdiff}) of the functional is in force
such that, for every $g\in {\cal G}$ and $t\in (0,\infty)$,
\begin{equation} \label{e.t.expT/twn}
     \sqrt{n}\,\bigl(\hspace{6\ei}T(P_{n,t,g})-T(P)\bigr) =
     t \hspace{6\ei}\langle \kappa|g\rangle + \Lo(n^0)
\end{equation}
Therefore, the asy.~hypotheses $J$, $H$ and~$K$ concern
  $(g,t)\in {\cal G}\times(0,\infty)$ and may be expressed by
\begin{equation} \label{E:T:PinK/cone}
   J^0: g= 0\weg,\hspace{1.5em}
   J: \langle \kappa|g\rangle = 0\weg,\hspace{1.5em}
   H: \langle \kappa|g\rangle \le0 \weg,\hspace{1.5em}
   K: t \hspace{6\ei}\langle \kappa|g\rangle \ge c
\end{equation}
Depending on whether the tangent set~${\cal G}$ is a convex
cone~$\tilde{{\cal G}}$ or a linear space~$\bar{{\cal G}}$,
the hypotheses~$J$, $H$, and~$K$ will be denoted by $\tilde{J}$,
  $\tilde{H}$, and~$\tilde{K}$, respectively by~$\bar{J}$, $\bar{H}$,
and~$\bar{K}$; obviously, $\tilde{J}^0=\bar{J}^0=J^0$.
\par
Overparametrization
  $P_{n,t,g}=P_{n,t/\!\gamma,\gamma g}$ with $t,\gamma>0$,
  for $g\in {\cal G}$\,(a cone), is allowed in~(\ref{e.t.Qpath})
  but, in view of~(\ref{e.t.expT/twn}) and~(\ref{E:T:PinK/cone}),
consistent with the functional and the hypotheses $J^0$, $J$, $H$,~$K$.
Distinction of three (actually, five) null hypotheses $J^0$, $J$,
and~$H$ is essential to Theorem~\ref{t:t:1s-power} and
Proposition~\ref{p.t.sizextH1}.
\subsection{Asymptotic Power Bounds for Cones and Spaces}
Let us fix some level $\alpha\in (0,1)$, and denote by~$u_{\alpha}$ the
upper $\alpha$-point of the standard normal distribution function~$\Phi$,
such that $\Phi(-u_{\alpha})=\alpha$. 
We shall employ asy.~tests, that is, sequences of tests~$\tau_n$ at
sample size~$n$. Power and size of the tests~$\tau_n$ are going to be
evaluated under the $n$-fold product measures~$Q_n^n$ asy.,
as $n\to \infty$.
An asy.~test $(\hat\tau_n)$ is said to achieve an upper bound
 $ \Ninf_{K} \Nlimsup_{n}\itg \tau_n\,dQ_n^n \le \beta $
 with $\limsup_n$ replaced by $\liminf_n$, if itself fulfills the
side conditions on the test sequences~$(\tau_n)$ under consideration
and $ \Ninf_{K} \Nliminf_{n}\itg \hat\tau_n\,dQ_n^n = \beta $ holds.
\begin{Thm}\sl  \label{t:t:1s-power}
Let~$(\tau_n)$ be an asy.~test that maintains asy.~level\/~{\rm $\alpha$}
under\/~{\rm $J^0$,}
\begin{equation} \label{e:t:level}
  \Nlimsup_{n}\itg \tau_n \,dP^n\le \alpha \end{equation}
\begin{ABC} \item    \label{i.t.cone}
Then, in the case of a convex tangent cone\/~{\rm $\tilde{{\cal G}}$,}
\begin{equation} \label{e.t.pow.cone}
  \Ninf_{\tilde{K}} \Nlimsup_{n}\itg \tau_n\,dQ_n^n \le
  \Phi \Bigl(-u_{\alpha}+\frac{c}{\Vert \tilde{\kappa}\Vert} \:\Bigr)
\end{equation}
The upper bound\/~{\rm (\ref{e.t.pow.cone}),}
with $\limsup_n$ replaced by\/~{\rm $\liminf_n$,}
is achieved by the asy.~test 
\begin{equation}\label{E:T:topt/cone}
   \tilde{\tau}_n=\hJc \bigl (\rave_{i=1}^n
   \tilde{\kappa}(x_i) > \Vert \tilde{\kappa}\Vert
   \hspace{6\ei} u_{\alpha}\bigr )
\end{equation}
\item  \label{i.t.lin}
In the case of a linear tangent space\/~{\rm $\bar{{\cal G}}$,}
\begin{equation} \label{e.t.pow.lin}
  \Ninf_{\bar{K}} \Nlimsup_{n}\itg \tau_n\,dQ_n^n \le
  \Phi \Bigl(-u_{\alpha}+\frac{c}{\Vert \bar{\kappa}\Vert} \:\Bigr)
\end{equation}
The upper bound\/~{\rm (\ref{e.t.pow.lin}),}
with $\limsup_n$ replaced by\/~{\rm $\liminf_n$,}
is achieved by the asy.~test 
\begingroup \mathsurround0em\arraycolsep0em \begin{eqnarray}
\noalign{\vspace{-\abovedisplayskip}\vspace{\abovedisplayshortskip}}
\label{E:T:topt/lin}  & \Ds
   \bar{\tau}_n=\hJc \bigl ( \rave_{i=1}^n
   \bar{\kappa}(x_i) > \Vert \bar{\kappa}\Vert \hspace{6\ei}
    u_{\alpha}\bigr ) & \\
\noalign{\noindent Moreover,\nopagebreak}
\label{E:T:topt/lin:alfa} & \Ds
  \Nsup_{\bar{H}}\Nlimsup_{n}\itg \bar\tau_n\,dQ_n^n
      \le \alpha & \end{eqnarray}\endgroup \end{ABC}\end{Thm}
\begin{Bew} \enskip\mbox{} 
\begin{ABC}\item[{\rfi{\ref{i.t.cone}}}]
Given any $g\in \tilde{{\cal G}}$ such that $\langle \kappa|g\rangle>0$,
  put $ t_g = c\big/\!\langle \kappa|g\rangle $ and test $J^0$ vs.~the
simple subhypothesis~$(P_{n, t_g ,g}^n)$ of~$\tilde{K}$.
Path differentiabilty~(\ref{e:i:pathdiff}) ensures the following
well-known asy.~expansion of loglikelihoods under~$P^n$, 
\begin{equation} \label{E:T:llh}
     \log \frac{dP_{n, t_g ,g}^n}{dP^n}=
     t_g \,\rave_{1}^n g(x_i) -
    \Tfrac{1}{2}  t_g^2 \hspace{6\ei}{\Vert g\Vert}^2 + \Lo_{P^n}(n^0)
\end{equation}
Thus Corollary~3.4.2 of Rieder\footnote{HR, subsequently}(1994)
is in force and bounds the asy.~power under~$P_{n,t_g,g}^n$
subject to~(\ref{e:t:level}) from above by
  $ \Phi (-u_{\alpha} +t_g \hspace{6\ei}\Vert g\Vert\,) $.
Now let $g\in \tilde{{\cal G}}$ approach~$\tilde{\kappa}$ in~$L_2(P)$.
Then $ t_g \hspace{6\ei}\Vert g\Vert $ tends to
     $ c \hspace{6\ei} \Vert \tilde{\kappa}\Vert
       \big/\!\langle \kappa|\tilde{\kappa}\rangle =
       c \big/ \Vert \tilde{\kappa}\Vert $
where we have used that
  $ \langle \kappa|\tilde{\kappa}\rangle
    = \Vert \tilde{\kappa}\Vert^2 $,
and bound~(\ref{e.t.pow.cone}) is obtained as the limit
\begin{equation}
   \lim_{g\to \tilde{\kappa}}
   \Phi (-u_{\alpha} + t_g \hspace{6\ei}\Vert g\Vert\,)
   = \Phi \Bigl(-u_{\alpha} + \frac{c}{\Vert \tilde{\kappa}\Vert} \,\Bigr)
\end{equation}
\par      Towards achieving bound~(\ref{e.t.pow.cone}) by the
tests~$\tilde{\tau}_n$, the sums 
   $\rave_1^n\tilde{\kappa}(x_i)$ are, for every
   $(g,t)\in \tilde{{\cal G}}\times(0,\infty)$ asy.~normal
under~$P^n_{n,t,g}$,
\begin{equation}  
    \bigl(\rave_{1}^n \tilde{\kappa}(x_i) \bigr)
    ( P^n_{n,t,g} ) \gwto
    {\cal N}\bigl( \hspace{6\ei} t \langle \tilde{\kappa}|g \rangle,
    \Vert \tilde{\kappa}\Vert^2 \hspace{9\ei}\bigr)
\end{equation}
by~(\ref{E:T:llh}) and a LeCam lemma, confer HR~(1994; Corollary~2.2.6),
and so
\begin{equation} \label{E:T:aspow:tfi}
   \Nlim_{n}\itg \tilde{\tau}_n\, dP^n_{n,t,g} =
   \Phi \Bigl( -u_{\alpha}+ \frac{t \langle
   \tilde{\kappa}|g \rangle}{\Vert \tilde{\kappa}\Vert }\,\Bigr)
\end{equation}
Under~$J^0:g=0$, this limit equals~$\alpha$.
If $(P^n_{n,t,g})\in \tilde{K}$ then, by~(\ref{E:T:PinK/cone}), as $t>0$,
and since $ \langle\tilde{\kappa}|g\rangle \ge \langle \kappa|g\rangle
            \enskip \forall g\in \tilde{{\cal G}}$, also
  $ t \hspace{6\ei}\langle \tilde{\kappa}|g\rangle \ge
    t \hspace{6\ei}\langle \kappa|g\rangle \ge c $. Hence
\begin{equation}
   \Ninf_{\tilde{K}}\Nlim_{n}\itg \tilde{\tau}_n\, dP^n_{n,t,g}  \ge
   \Phi \Bigl( -u_{\alpha} + \frac{c}{\Vert \tilde{\kappa}\Vert }\,\Bigr)
\end{equation}
\item[{\rfi{\ref{i.t.lin}}}]
With $\bar{\kappa}$ and~$\bar{K}$ in the place of~$\tilde{\kappa}$
and~$\tilde{K}$, the proof of bound~(\ref{e.t.pow.lin}) is the same as
in case~(\ref{i.t.cone}). The limit corresponding to~(\ref{E:T:aspow:tfi})
for the tests~$\bar{\tau}_n$ is 
\begin{equation} \label{E:T:limpow/g}
   \Nlim_{n}\itg \bar{\tau}_n\, dP^n_{n,t,g} =
   \Phi \Bigl( -u_{\alpha} +\frac{t \langle
   \bar{\kappa}|g \rangle}{\Vert \bar{\kappa}\Vert }\,\Bigr) =
   \Phi \Bigl( -u_{\alpha}+\frac{t \langle
   \kappa|g \rangle}{\Vert \bar{\kappa}\Vert }\,\Bigr)
\end{equation}
since $ \kappa- \bar{\kappa}\perp \bar{{\cal G}}$.
If $(P^n_{n,t,g})\in \bar{H}$, then $ t \langle \kappa|g\rangle \le 0 $
by~(\ref{E:T:PinK/cone}) and $t>0$. Therefore
\begin{equation}
   \Nlim_{n} \itg \bar{\tau}_n\, dP_{n,t,g}^n
   = \Phi \Bigl( -u_{\alpha}+ \frac{t \langle
   \kappa|g \rangle}{\Vert \bar{\kappa}\Vert }\,\Bigr) \le
     \Phi ( -u_{\alpha}+0) = \alpha
\end{equation}
is obtained from~(\ref{E:T:limpow/g}), and
proves~(\ref{E:T:topt/lin:alfa}).\qed                 \end{ABC} \end{Bew}
\begin{Rem}\rm  Although Theorem~\ref{t:t:1s-power}\,(\ref{i.t.cone}),
for convex tangent cones, is straightforward to prove, it seems to have
been omitted in literature so far. In its proof, 
  $\tilde{\kappa}$~acts as a limiting least favorable tangent,
  as does~$\bar{\kappa}$ in the proof            
of Theorem~\ref{t:t:1s-power}\,(\ref{i.t.lin}).
The latter result, for linear tangent spaces, compares with
Pfanzagl and Wefelmeyer (1982; chapter~8),
van der Vaart~(1998; Theorem~25.44, Lemma~25.45), as well as
Beran~(1983; Theorem~1) and HR~(1994; Theorem~4.3.8) who,
in robust statistics, encounter linear tangent spaces~$\bar{{\cal G}}$
with maximal closure
  $\cl\bar{{\cal G}}=L_2(P)\cap \nolinebreak\{1\}^\perp$
(so that $\bar{\kappa}= \kappa-\Ew \kappa$ there).\qed
\end{Rem} 
\subsection{Comparison of Cones and Their Linear Spans}
                                                \label{ss.t.pow.cone/lin}
Let us consider~$P$ a member of two models
  $\tilde{{\cal P}}\subset\bar{{\cal P}}$ whose tangent sets at~$P$
are a convex cone~$\tilde{{\cal G}}$, respectively the linear span
of~$\tilde{{\cal G}}$,
\begin{equation} \label{e.t.Gbar=spanGtil}
  \bar{{\cal G}}=\lin \tilde{{\cal G}} 
\hspace{-1.25em}\end{equation}
\paragraph{Power Comparison}
In this situation, we have $\tilde{J}\subset \bar{J}$,
  $\tilde{H}\subset \bar{H}$, $\tilde{K}\subset \bar{K}$,
and~$J^0$ should be easier to test 
vs.~$\tilde{K}$ than vs.~$\bar{K}$. In fact,
{\samepage
\begingroup \mathsurround0em\arraycolsep0em \begin{eqnarray}
\label{e.t.powtk>bk}
   \Phi \Bigl(-u_{\alpha}+\frac{c}{\Vert \tilde{\kappa}\Vert} \,\Bigr)
   &{}>{}&
   \Phi \Bigl(-u_{\alpha}+\frac{c}{\Vert \bar{\kappa}\Vert} \,\Bigr) \\
\noalign{\vspace{-.5ex}\noindent because \nopagebreak \vspace{-.5ex}}
\label{e.t.normtk<bk}
   \Vert \tilde{\kappa}\Vert & {}< {}& \Vert \bar{\kappa}\Vert
\end{eqnarray}\endgroup
unless $ \bar{\kappa} \in \cl\tilde{{\cal G}}$,
in which case $ \tilde{\kappa}=\bar{\kappa} $ and
the two power bounds coincide. }
\par
This is a consequence of
  $ \Vert \tilde{\kappa}\Vert^2 = \langle \kappa|\tilde{\kappa} \rangle
    = \langle \bar{\kappa}|\tilde{\kappa} \rangle $
and the Cauchy--Schwarz inequality:
  $  \langle \bar{\kappa}|\tilde{\kappa} \rangle \le
     {\Vert \bar{\kappa}\Vert} \hspace{6\ei}\Vert \tilde{\kappa}\Vert $,
where equality holds iff~$\tilde{\kappa}$ is some positive multiple
of~$\bar{\kappa}$, in which case
$\bar{\kappa}\in \cl\tilde{{\cal G}}$ and $\tilde{\kappa}=\bar{\kappa}$.
\paragraph{Sample Size Comparison}
Allowing for different sample sizes~$\tilde{n}$ and $\bar{n}$,
respectively, such that $\tilde{n}/n\to \tilde{\gamma}$ and
       $\bar{n}/n\to \bar{\gamma}$ for some
$\tilde{\gamma}, \bar{\gamma}\in (0,\infty)$, the asy.~power bounds
(\ref{e.t.pow.cone}) and~(\ref{e.t.pow.lin}) are the same iff
\begin{equation} \label{e.t.nt/nb}
      \bar{\gamma}:\tilde{\gamma} =
      {\Vert \bar{\kappa}\Vert}^2 \!: {\Vert \tilde{\kappa}\Vert}^2
\end{equation}
Thus, observations at the higher rate
  $ \Vert \bar{\kappa}\Vert^2\!\big/\Vert \tilde{\kappa}\Vert^2 $ are needed
by~$(\bar{\tau}_{\bar{n}})$ to achieve, subject to level~$\alpha$ on~$J^0$,
the same power vs.~$\bar{K}$ as~$(\tilde{\tau}_{\tilde{n}})$ vs.~$\tilde{K}$.
\begin{Exa}\rm  \label{ex1.t.kbkt/cone}
Consider the standard normal $P= {\cal N}(0,1)$ and $\kappa(x)=x$
the identity on the real line; $\kappa$~is 
the influence curve at~$P$ of the expectation functional
as well as of the one-sample normal scores rank functional,
\begingroup \mathsurround0em\arraycolsep0em \begin{eqnarray}
   & \Ds  E(Q) = \itg_{-\infty}^{\infty} x\,Q(dx)  & \\ \hspace{-2.25em}
   & \Ds  R(Q) = 2\itg_{0}^{\infty} \Phi^{-1} \bigl( \Tfrac{1}{2}
               + \Tfrac{1}{2} [ Q(x)-Q(-x) ] \bigr)\,Q(dx)
               - 2 \hspace{9\ei}\varphi(0)   &
\end{eqnarray}\endgroup
where $\varphi= \nolinebreak \dot \Phi$ denotes the standard normal density,
and~$Q(x)=Q \bigl((-\infty,x]\hspace{6\ei}\bigr)$.
\par
As tangents at~$P$, consider the sign-function $g_1(x)=\zi(x)$
and the function $g_2(x)=\mu \hspace{3\ei}\zi(x)\Jc(|x|\le a)$ with
$\mu, a\in (0,\infty)$. Then $\Vert g_1\Vert=1=\Vert \kappa\Vert $,
and~$\mu=\mu_a$ may be determined by $\mu^{-2}_a=2 \hspace{6\ei}\Phi(a)-1$
such that also $\Vert g_2\Vert=1$. Then the coefficients
  $b_i=\nolinebreak \langle \kappa|g_i \rangle $
  and $c= \langle g_1|g_2\rangle$ are given by
\begin{equation} \label{e.t.ex1.coeff}
  b_1=2 \hspace{6\ei}\varphi(0)\weg,\hspace{1.5em}  b_2=2 \hspace{6\ei}
  \mu \bigl[\varphi(0)-\varphi(a)\bigr]\weg,\hspace{1.5em}
  c= 2 \hspace{6\ei}\mu \bigl[ \Phi(a)-\Tfrac{1}{2}\hspace{3\ei}\bigr]
\end{equation}
As tangent sets at~$P$, employ the (closed) convex
cone~$\tilde{{\cal G}}=\cl \tilde{{\cal G}}$ and (closed) linear
space~$\bar{{\cal G}}=\cl \bar{{\cal G}}=\lin \tilde{{\cal G}}$ spanned
by the tangents $g_1$ and~$g_2$,
\begin{equation}
  \tilde{{\cal G}}= \bigl\{\, \gamma_1 g_1 + \gamma_2 \hspace{3\ei}g_2
                    \bigm|  \gamma_i\ge0 \,\bigr\}  \weg,\qquad
  \bar{{\cal G}}=   \bigl\{\, \gamma_1 g_1 + \gamma_2 \hspace{3\ei}g_2
                    \bigm| \gamma_i\in\R\,\bigr\}
\end{equation}
Via~(\ref{e.i.dPs1fach}), the cone~$\tilde{{\cal G}}$ defines a set of
positively asymmetric alternatives to~$P$.
\par     Unconstrained minimization of
  $\Vert{\kappa-\gamma_1\hspace{2\ei}g_1-\gamma_2\hspace{4\ei}g_2}\Vert$
being equivalent to the orthogonality relations
  $ \gamma_1+\gamma_2 \hspace{3\ei}c = b_1$ and
  $ \gamma_1 c + \gamma_2 =b_2$, the canonical gradient is
\begin{equation} \label{e.t.ex1.kb}
  \bar{\kappa}= \bar{\gamma}_1 g_1+\bar{\gamma}_2 \hspace{3\ei}g_2
  \qquad \mbox{where}\quad
   \bar{\gamma}_1= \frac{b_1-b_2 \hspace{6\ei}c}{1-c^2}\weg,\enskip
    \bar{\gamma}_2= \frac{b_2-b_1c}{1-c^2}
\end{equation}
In the appendix we show that $\bar{\gamma}_1>0>\bar{\gamma}_2$;
hence $\bar{\kappa}\in \bar{{\cal G}}\setminus \tilde{{\cal G}}$.
\par
The constrained minimization of
  $\Vert{\kappa-\gamma_1g_1-\gamma_2\hspace{3\ei}g_2}\Vert$
subject to $\gamma_i\ge0$ is a convex and well-posed problem;
HR~(1994; Theorem~B.2.3, Definition~B.2.9). Thus there exist
multipliers $\beta_i\ge0$ such that the solutions $\tilde{\gamma}_i\ge0$
minimize the following Lagrangian over $\gamma_i\in\R$,
\begin{equation} \label{e.t.ex1.Lagr}
\mathsurround0em\arraycolsep0em \begin{array}{l}
\Ds {\Vert{\kappa-\gamma_1g_1-\gamma_2\hspace{3\ei}g_2}\Vert}^2
    - 2 \hspace{6\ei}\beta_1 \gamma_1
    - 2 \hspace{6\ei}\beta_2 \hspace{3\ei}\gamma_2 - \mbox{const}\\
\Ds \rule{0pt}{3ex}\hspace{2.25em} {}=
   \bigl[\gamma_1-(b_1+\beta_1)\bigr]^2 +
   \bigl[\gamma_2-(b_2+\beta_2)\bigr]^2 +
    2 \hspace{6\ei}c \hspace{12\ei}\gamma_1 \gamma_2
\end{array}\end{equation}
Moreover, $\beta_i \tilde{\gamma}_i=0$.
Since $\tilde{\kappa}\ne \bar{\kappa}$,
not both~$\beta_0$ and~$\beta_1$ can vanish.
\par
In case $\beta_1>0$ we obtain that $\tilde{\gamma}_1=0$ and
  $\tilde{\gamma}_2=b_2+\beta_2$, where $\beta_2=0$ because
  $\beta_2 \tilde{\gamma}_2=0$ and $b_2\ge0$. Hence $\tilde{\gamma}_2=b_2$
and $\Vert \kappa-b_2 \hspace{3\ei}g_2\Vert^2=1-b_2^2  $.
Likewise, if $\beta_2>\nolinebreak 0$ we obtain that $\tilde{\gamma}_2=0$
and $\tilde{\gamma}_1=b_1+\beta_1$, where $\beta_1=0$ because
  $\beta_1 \tilde{\gamma}_1=0$, hence $\tilde{\gamma}_1=b_1$
  and $\Vert \kappa-b_1g_1\Vert^2=1-b_1^2  $. Since $b_2<b_1$,
we have thus proved that $\tilde{\kappa}=b_1g_1$ always.
\par
Numerical values for $a=1$ are 
\begin{equation} \mathsurround0em\arraycolsep0em \begin{array}{c}
\Ds \mu=1.210\weg,\hspace{36\ei}
    b_1=0.798\weg,\hspace{36\ei}
    b_2=0.380 \weg,\hspace{36\ei}
    c=0.826 \\ 
\Ds \rule{0pt}{3ex}
    \bar{\gamma}_1=1.525\weg,\hspace{36\ei}
    \bar{\gamma}_2=-0.880 \weg,\hspace{36\ei} 
    {\Vert \bar{\kappa}\Vert}^2 = 0.882\weg,\hspace{36\ei}
    {\Vert \tilde{\kappa}\Vert}^2 = 0.637 \\ 
\Ds \rule{0pt}{3ex}
    \Vert \bar{\kappa}\Vert^2 \!: \Vert \tilde{\kappa}\Vert^2
    = 1.386\weg,\hspace{36\ei}
    \Vert \tilde{\kappa}\Vert^2 \!:\Vert \bar{\kappa}\Vert^2 = .721
\end{array}\end{equation}
The value~$.721$, to the third digit, turns out to be the minimum of
  $ \Vert \tilde{\kappa}\Vert^2\! \big/ \Vert \bar{\kappa}\Vert^2 $
with respect to~$a\in (0,\infty)$.%
\qed\end{Exa}  
\subsection{\boldmath
            Level Breakdown of~$(\tilde{\tau}_n)$}    \label{ss.t.Lbruch}
In the setup~(\ref{e.t.Gbar=spanGtil}):
   $\bar{{\cal G}}=\lin \tilde{{\cal G}}$, in view of~(\ref{E:T:topt/lin:alfa}),
the tests~$\bar{\tau}_n$ automatically maintain asy.~level~$\alpha$ on the
left-sided extension~$\bar{H}$ of~$\bar{J}$ and~$J^0$,
where $\bar{H} \supset \tilde{H}$.
On the contrary, the analogue to~(\ref{E:T:topt/lin:alfa}) for extensions
  $\tilde{H}\supset \tilde{J}$ of~$J^0$ and the tests~$\tilde{\tau}_n$ can
in general not be achieved.
\par
Note that
\begin{equation} \label{e.t.kbnekt.sprod.g}
   \bar{\kappa}\ne \tilde{\kappa}\iff \exists\,g\in \tilde{{\cal G}}:
   \langle \kappa|g \rangle < \langle \tilde{\kappa}|g \rangle
\end{equation}
\begin{Prop}\sl  \label{p.t.sizextH1}
Assume the convex cone~$\tilde{{\cal G}}$ contains a tangent~$g_0$ such that
\begingroup \mathsurround0em\arraycolsep0em \begin{eqnarray}
\label{e.t.kg<0<tlkg} & \Ds \langle \kappa|g_0 \rangle \le 0
                            < \langle \tilde{\kappa}|g_0 \rangle & \\
\noalign{\noindent Then\nopagebreak} & \Ds
  \Nsup_{\tilde{J}}\Nlimsup_{n}\itg \tilde\tau_n\,dQ_n^n = 1 &
\end{eqnarray}\endgroup \end{Prop}
\begin{Bew} \enskip
If $\langle \kappa|g_0 \rangle = 0 $, then
      $(P^n_{n,t,g_0})\in\tilde{J}\enskip \forall t\in (0,\infty)$.
In view of~(\ref{E:T:aspow:tfi}), therefore, the tests $\tilde{\tau}_n$
have asy.~size at least
\begin{equation}
     \Nsup_{t>0} \Nlim_{n}\itg \tilde{\tau}_n\, dP^n_{n,t,g_0} =
     \sup_{t>0} \: \Phi \Bigl( - u_{\alpha} + \frac{t \langle
     \tilde{\kappa}|g_0 \rangle}{\Vert \tilde{\kappa}\Vert }\,\Bigr)\, = 1
\end{equation}
because
  $ \lim_{t\to \infty}t\langle \tilde{\kappa}|g_0 \rangle = \infty $
due to $ \langle \tilde{\kappa}|g_0 \rangle > 0$.
\par     In case
  $\langle \kappa|g_0 \rangle < 0 < \langle \tilde{\kappa}|g_0 \rangle $,
a suitable convex combination~$g_{01}$ of~$g_0$ and~$\tilde{\kappa}$, since
  $ 0<\langle \kappa|\tilde{\kappa}\rangle = \Vert \tilde{\kappa}\Vert^2 $,
will satisfy $  \langle \kappa|g_{01} \rangle = 0 <
                \langle \tilde{\kappa}|g_{01} \rangle  $.\qed \end{Bew}
\begin{Exa}\rm   \label{ex2.t.JnotH/cone}
In Example~\ref{ex1.t.kbkt/cone}, although $\bar{\kappa}\ne \tilde{\kappa}$,
condition~(\ref{e.t.kg<0<tlkg}) is not fulfilled, because $b_1,b_2>0$,
and so $ \langle \kappa|g \rangle \le 0 $ can hold for
       $g\in \tilde{{\cal G}}$ only if $g=0$.
\par     However, in the setup of Example~\ref{ex1.t.kbkt/cone},
have tangent~$g_2$ be replaced by the function
\begin{equation} \label{e.t.ex2.g3}
   g_3(x) = - g_3(-x) =
          \cases{ \delta  & \mbox{if $0<x\le a$}\cr
                 - \eta & \mbox{if $a<x$}}
\end{equation} 
with $a,\delta, \eta\in (0,\infty)$.
In the appendix we show that, given any $a\in (0,\infty)$, the constants
  $\eta=\eta_a$ and $ \delta = \delta_a $ may be determined by
  $ \delta_a = \sigma_a \hspace{6\ei} \eta_a $ and
\begin{equation} \label{e.t.sada}
  \eta_a^{-2} = 2 \hspace{6\ei}
     \bigl( \sigma_a^2 \bigl[ \Phi(a) - \Tfrac{1}{2}\hspace{3\ei}\bigr]
     + \bigl[ 1-\Phi(a)\bigr]\bigr) 
  \weg, \qquad 
  \sigma_a = a \hspace{9\ei} \frac{1-\Phi(a)}{\varphi(0)-\varphi(a)}
\end{equation}
Then $\Vert g_3\Vert =1$ and
\begin{equation}
   \langle \kappa|g_3 \rangle <0< \langle g_1|g_3 \rangle
\end{equation}
By the method of Lagrange multipliers, in the appendix, we prove that
\begin{equation}
   \tilde{\kappa}= \langle \kappa|g_1 \rangle \hspace{6\ei}g_1
\end{equation}
where $ \langle \kappa|g_1 \rangle = 2 \hspace{6\ei}\varphi(0)>0$,
and so $ \langle \kappa|g_3 \rangle < 0 <
          \langle \kappa|g_1 \rangle \langle g_1|g_3 \rangle =
            \langle \tilde{\kappa}|g_3 \rangle $,
which implies~(\ref{e.t.kg<0<tlkg}) for $g_0=g_3$.
\qed\end{Exa}
Making use of the following uniqueness result (Proposition~\ref{p.t.u}),
we conclude that testing the slightly bigger null
hypothesis~$\tilde{J}\supset \nolinebreak J^0$, or the even bigger
one-sided extension~$\tilde{H}$ of~$\tilde{J}$, vs.~$\tilde{K}$,
is inevitably bound to larger error probabilities than those
given in Theorem~\ref{t:t:1s-power}(\ref{i.t.cone}) for testing
$J^0$ vs.~$\tilde{K}$. This is contrary to the extension of~$J^0$
to~$\bar{J}$ and~$\bar{H}$, vs.~$\bar{K}$, which goes for free
in Theorem~\ref{t:t:1s-power}(\ref{i.t.lin}).
\begin{Rem}\rm    \label{r.t.extH.test/cone}
The minimum asy.~power
  $ \Phi \bigl(-u_{\alpha} + c \big/\Vert \bar{\kappa}\Vert \,\bigr) $
achieved by the asy.\ test~$(\bar{\tau}_n)$ under~$\bar{K}$ stays
the same under~$\tilde{K} \subset \nolinebreak\bar{K}$, that is,
does not increase,
\begin{equation}
  \Ninf_{\tilde{K}} \Nlim_{n}\itg \bar{\tau}_n\,dQ_n^n =
  \Phi \Bigl(-u_{\alpha}+\frac{c}{\Vert \bar{\kappa}\Vert} \:\Bigr)
\end{equation}
Indeed, pick any $g\in \tilde{{\cal G}}$ such that
  $\langle \kappa|g\rangle>0$; for example, $g=\tilde{\kappa}$ itself.
Then choose $t\in (0,\infty)$ such that
  $ t \hspace{6\ei}\langle \kappa|g\rangle = c $,
and apply (\ref{E:T:PinK/cone}) and~(\ref{E:T:limpow/g}).
\par     Whether
  $ \Phi \bigl(-u_{\alpha} + c \big/\Vert \bar{\kappa}\Vert \,\bigr) $ is
the largest minimum asy.~power that can be achieved vs.~$\tilde{K}$,
subject to asy.~level~$\alpha$ under~$\tilde{H}$, respectively
only under~$\tilde{J}$, is unknown.
In particular, we do not know if there exists some
function~$ \eta\in \nolinebreak L_2(P)$ of smaller norm 
  $ \Vert \eta\Vert < \Vert \bar{\kappa}\Vert $ and such that,
  for each $ g\in \tilde{{\cal G}} $,
\begin{equation} \label{e.t.sprodK<0>}
\mathsurround0em \arraycolsep0em 
  \left. \begin{array}{rrcl}
  \Ds \tilde{J}:{} & \Ds \langle \kappa|g\rangle &{}={}& 0 \\
  \Ds \tilde{H}:{} & \Ds \langle \kappa|g\rangle &{}\le{}& 0
          \end{array}\:\right\}
          \so \langle \eta|g\rangle \le0 \weg, \qquad
   \langle \kappa|g\rangle > 0 \so
   \langle \eta|g\rangle \ge \langle \kappa|g\rangle
\end{equation}
In connection with asy.~median unbiased, two-sided confidence
limits for cones, the corresponding function~$\eta$ cannot exist;
confer Subsection~\ref{ss.e.vs}, where instead of~(\ref{e.t.sprodK<0>})
the simpler condition~(\ref{e.e.sprodk.1/clincone}) occurs.
\qed\end{Rem}
\subsection{Uniqueness of Most Powerful Tests}
In the setup of Theorem~\ref{t:t:1s-power}, the optimal tests
  $\tilde{\tau}_n$ and $\bar{\tau}_n$ defined by~(\ref{E:T:topt/cone})
and~(\ref{E:T:topt/lin}), respectively, are unique up to
terms~$\Lo_{P^n}(n^0)$ tending stochastically to zero under~$P^n$.
\begin{Prop}\sl \label{p.t.u}  
Suppose that an asy.~test~$(\tau_n)$ satisfies\/~{\rm (\ref{e:t:level}),}
and achieves the asy.~power bound\/~{\rm (\ref{e.t.pow.cone})}
in case\/~{\rm (\ref{i.t.cone}),} respectively
bound\/~{\rm (\ref{e.t.pow.lin})} in case\/~{\rm (\ref{i.t.lin})}.
Then necessarily 
\begin{equation} \label{e.t.uopt}
   \tau_n= \cases{\Ds \tilde{\tau}_n + \Lo_{P^n}(n^0)
           & in case\/~{\rm (\ref{i.t.cone}),} \enskip respectively\cr
           \Ds \bar{\tau}_n + \Lo_{P^n}(n^0)
           & in case\/~{\rm (\ref{i.t.lin}).} }
\end{equation}
Conversely, form\/~{\rm (\ref{e.t.uopt})} implies that the
asy.\ test\/~{\rm $(\tau_n)$} satisfies\/~{\rm (\ref{e:t:level})}
and achieves bound\/~{\rm (\ref{e.t.pow.cone}),} respectively
satisfies\/~{\rm (\ref{E:T:topt/lin:alfa})} and
achieves bound\/~{\rm (\ref{e.t.pow.lin}).}
\end{Prop}
\begin{Bew}\enskip
In regard of the proof to Theorem~\ref{t:t:1s-power}, the proposition
is a straightforward consequence of the uniqueness result in
HR~(1994; Corollary~3.4.2\,b with $\sigma>0$)---under the provision
however that $\tilde{\kappa}\in \tilde{{\cal G}}$ and
     $\bar{\kappa}\in \bar{{\cal G}}$, respectively.
Since, in general, the tangent set~${\cal G}$ (convex cone or linear space)
needs not be closed, we have to incorporate an approximation
in~$L_2(P)$ of~$\tilde{\kappa}$ and~$\bar{\kappa}$ by elements
of~$\tilde{{\cal G}}$ and~$\bar{{\cal G}}$, respectively.
This is the reason for the following proof.
\par     \smallskip   
Thus, given any $t\in (0,\infty)$ and $h\in L_2(P)$, $h\ne0$, $\Ew h=0$,
we shall show that~(\ref{e:t:level}) and
\begingroup \mathsurround0em\arraycolsep0em \begin{eqnarray}
\noalign{\vspace{-\abovedisplayskip}\vspace{\abovedisplayshortskip}}
\label{e.t.approxopt} & \Ds
  \liminf_{(s,g)\to (t,h)} \Nliminf_{n} \itg \tau_n\,dP^n_{n,s,g}
  \ge \Phi \bigl(-u_{\alpha}+t \hspace{3\ei \Vert h\Vert }\,\bigr) & \\
\noalign{\noindent imply that \nopagebreak} \label{e.t.1exp} & \Ds
   \tau_n=\hJc \bigl(\rave_{1}^n h(x_i) >
   \Vert h\Vert \hspace{6\ei} u_{\alpha}\bigr)
   + \Lo_{P^n}(n^0) & \end{eqnarray}\endgroup
In proving this, it is no restriction to set $s=t=1$,
and then delete $s$ and~$t$ from notation;
in particular, we write $P_{n,s,g}=P_{n,g}$.
\par     Given any $\delta\in (0,1)$, $\delta< \Vert h\Vert $, choose~$g$
so close to~$h$ that
\begin{equation} \label{e.t.nbl<d}
     {\Vert g-h\Vert}^2 < \delta^3\weg,\hspace{.33em}
     \bigl|\Vert g\Vert^2- \Vert h\Vert^2  \bigr|<
     2 \hspace{6\ei}\delta 
     \hspace{.75em}\mbox{and}\hspace{.75em}
     |\hspace{9\ei}\beta_g- \beta_h|<\delta\weg,\hspace{.33em}
     |\hspace{9\ei}\ell_g- \ell_h|<\delta
\end{equation}
for the norm based quantities
    $\beta_g=\Phi \bigl(-u_{\alpha}+ \Vert g\Vert\,\bigr)$ and
    $ \ell_g=\Vert g\Vert \hspace{4\ei} u_{\alpha}
                    - \frac{1}{2} \Vert g\Vert^2 $,
and such that, making use of~(\ref{e.t.approxopt}), moreover
\begin{equation} \label{e.t.powapprox}
     \Nliminf_{n}\itg \tau_n\,dP^n_{n,g} \ge \beta_h - \delta
\end{equation}
The proof employs the following Neyman--Pearson tests~$\tau_{n,g}^*$ for
   $P^n$ vs.~$P^n_{n,g}$,
\begin{equation}
    \tau_{n,g}^*= \Jc (L_{n,g}>\ell_g) \weg,\qquad
    L_{n,g}=\log dP^n_{n,g}/dP^n
\end{equation}
As the loglikelihoods $L_{n,g}$ are asy.\
   $ {\cal N}\bigl(- \frac{1}{2}\Vert g\Vert^2,
      \Vert g\Vert^2 \hspace{6\ei} \bigr) $ under~$P^n$, 
\begin{equation} \label{e.t.NPTg.abg}
    \alpha_n=\itg \tau_{n,g}^* \,dP^n \longrightarrow \alpha\weg,
    \qquad
    \beta_n=\itg \tau_{n,g}^* \, dP_{n,g}^n \longrightarrow \beta_g
\end{equation}
By (\ref{e:t:level}), (\ref{e.t.nbl<d}), and~(\ref{e.t.powapprox}),
some $n_0=n_0(\delta)$ exists such that for all $n\ge n_0$,
\begin{equation} \label{e.t.taun.anbndelta}
  \itg \tau_n\, dP^n \le \alpha_n+ 3 \hspace{6\ei}\delta\weg,\qquad
  \itg \tau_n\,dP^n_{n,g} \ge \beta_n - 3 \hspace{6\ei}\delta
\end{equation}
Then Lemma~\ref{l.x.u} tells us that, for all such $n\ge n_0$
and for every $\varepsilon\in (0,1)$,
\begingroup \mathsurround0em\arraycolsep0em
\begin{eqnarray}
 & \Ds  |\nu_{n,g}| \bigl\{|\tau_n- \tau_{n,g}^*|>\varepsilon\bigr\} \le
 3\hspace{6\ei}(1+c_{g})\hspace{1\ei}\Tfrac{\Ds\delta}{\Ds\varepsilon} & \\
\noalign{\noindent where \nopagebreak}
 & \Ds  \nu_{n,g}= P^n_{n,g}-c_g P^n\weg,\qquad c_g=e^{\ell_g}     &
\end{eqnarray}\endgroup
Fix $\varepsilon\in (0,1)$ and set
    $ A_{n,g} =\bigl\{|\tau_n- \tau_{n,g}^*|>\varepsilon\bigr\}$.
Fix any $\rho\in (0,1)$. Then the probability
    $ P^n \bigl( A_{n,g}\cap \{L_{n,g}>\ell_g+\rho\}\bigr) $
is bounded above by
\begingroup \mathsurround0em\arraycolsep0em
\begin{eqnarray}
\hspace{-2.25em}  \frac{1}{e^{(\ell_g+\rho)}-e^{\ell_g}}
  \int_{A_{n,g}\cap \{L_{n,g}>\ell_g+\rho\}}
  \bigl( e^{L_{n,g}}- e^{\ell_g}\bigr) \,dP^n   & {}\le {} &
  \frac{|\nu_{n,g}|(A_{n,g})}{c_g \hspace{4\ei}(e^\rho-1)}     \\
\noalign{\vspace{\belowdisplayshortskip}\pagebreak[0]
\mathsurround=\msu\noindent   Likewise,
  $ P^n \bigl( A_{n,g}\cap \{L_{n,g}<\ell_g-\rho\}\bigr) $
  is bounded above by\nopagebreak\vspace{\abovedisplayskip}}
\hspace{-2.25em}  \frac{1}{e^{\ell_g}-e^{(\ell_g-\rho)}}
  \int_{A_{n,g}\cap \{L_{n,g}<\ell_g-\rho\}}
  \bigl( e^{\ell_g}-e^{L_{n,g}}\bigr) \,dP^n    & {}\le {} &
  \frac{|\nu_{n,g}|(A_{n,g})}{c_g \hspace{4\ei}(1-e^{-\rho})}
\end{eqnarray}\endgroup
Put $\eta_{\rho}=(e^\rho+1)\big/(e^\rho-1)$ and 
use $ |\hspace{9\ei}\ell_g- \ell_h|<\delta $, hence
    $ c_g>e^{-\delta}c_h$, 
to conclude that
\begin{equation} \mathsurround0em\arraycolsep0em
\hspace{-1em}\begin{array}{r@{{}\le{}}l} \Ds P^n(A_{n,g}) & \Ds
   3 \hspace{9\ei}\eta_{\rho} \hspace{6\ei} (1+c_g^{-1})
   \hspace{1\ei} \Tfrac{\Ds\delta}{\Ds\varepsilon} +
   P^n \bigl\{|L_{n,g}-\ell_g|\le \rho\bigr\} \\
\rule{0pt}{4ex} & \Ds
   3 \hspace{9\ei}\eta_{\rho} \hspace{6\ei} (1+e^{\delta}c_h^{-1})
   \hspace{1\ei} \Tfrac{\Ds\delta}{\Ds\varepsilon}
   + P^n \bigl\{|L_{n,g}-\ell_g|\le \rho\bigr\}
\end{array}\end{equation}
Asy.~normality of~$L_{n,g}$ under~$P^n$, and~(\ref{e.t.nbl<d})
ensuring $ \Vert g\Vert \ge \Vert h\Vert - \delta $, imply
\begin{equation}
   \Nlim_{n} P^n \bigl\{|L_{n,g}-\ell_g|\le \rho\bigr\}
   \le 2 \hspace{6\ei} \rho \hspace{6\ei} \frac{\varphi(0)}{\Vert g\Vert }
   \le 2 \hspace{6\ei} \rho \hspace{6\ei}
                       \frac{\varphi(0)}{\Vert h\Vert - \delta}
\end{equation}
It follows that, for all $\delta\in (0,1)$, $\delta<\Vert h\Vert $,
   and for all $\rho\in (0,1)$,
\begingroup \mathsurround0em\arraycolsep0em
\begin{eqnarray} \hspace{-2.25em} \label{e.t.PAasybound} & \Ds
   \limsup_{g\to h}\Nlimsup_{n} P^n (A_{n,g}) \le
   3 \hspace{9\ei}\eta_{\rho}\hspace{6\ei}
   (1+ e^{\delta} c_h^{-1}) \hspace{1\ei} \Tfrac{\Ds\delta}{\Ds\varepsilon}
   +  2 \hspace{6\ei} \rho \hspace{6\ei}
   \frac{\varphi(0)}{\Vert h\Vert - \delta}  & \\
\noalign{\noindent Hence \nopagebreak}
\label{e.t.gtohPAn0}   & \Ds    \lim_{g\to h}\Nlimsup_{n}
         P^n \bigl\{|\tau_n- \tau_{n,g}^*|>\varepsilon\bigr\} =0 &
\nopagebreak \end{eqnarray}\endgroup  \nopagebreak
if we first let $\delta$ and then $\rho$ approach~$0$
in~(\ref{e.t.PAasybound}).
\par     
Furthermore, comparing the Neyman--Pearson tests
         $\tau_{n,g}^*$ and~$\tau_{n,h}^*$, we get
\begin{equation} \label{e.t.tgzuth}
\mathsurround0em\arraycolsep0em \begin{array}{rcl}
\Ds  P^n \bigl\{ | \tau_{n,g}^*- \tau_{n,h}^* | > \varepsilon\bigr\}
     & {}\le{} & \Ds   P^n \bigl\{L_{n,g}>\ell_g\weg,\:
                       L_{n,h}<\ell_h -4 \hspace{6\ei} \delta\bigr\} \\
\rule{0pt}{3.3ex} & & \Ds \hspace{.75em} {} +
   P^n \bigl\{L_{n,g}\le \ell_g\weg,\:
                       L_{n,h}>\ell_h +4 \hspace{6\ei} \delta\bigr\} \\
\rule{0pt}{3.3ex} & & \Ds \hspace{2.25em} {} +
   P^n \bigl\{ | L_{n,h}-\ell_h|\le 4 \hspace{6\ei} \delta\bigr\}
\end{array}\end{equation}
The 3rd summand on the RHS, by the asy.~normality of~$L_{n,h}$
under~$P^n$, satisfies
\begin{equation}
  \Nlim_{n} P^n \bigl\{ | L_{n,h}-\ell_h|\le 4 \hspace{6\ei} \delta\bigr\}
  \le 8 \hspace{6\ei}\delta \hspace{9\ei} \frac{\varphi(0)}{\Vert h\Vert }
\end{equation}
The first two summands on the RHS in~(\ref{e.t.tgzuth}),
since $|\hspace{9\ei}\ell_g- \ell_h|<\delta$, are bounded by
  $ P^n \bigl\{ | L_{n,g}-L_{n,h}|> 3 \hspace{6\ei}\delta\bigr\} $.
Invoke the loglikelihood expansion~(\ref{E:T:llh}) and make use of
  $ \bigl|\Vert g\Vert^2- \Vert h\Vert^2  \bigr|< 2 \hspace{6\ei}\delta $
in order to bound
  $ P^n \bigl\{ | L_{n,g}-L_{n,h}|> 3 \hspace{6\ei}\delta\bigr\} $ by
\begin{equation}
  P^n \bigl\{\bigl| \rave_{1}^n {(g-h)(x_i)}\bigr|
  > 2 \hspace{6\ei}\delta  - \Lo_{P^n}(n^0)\biggr\}
\end{equation}
which, in turn, is bounded by some~$\Lo(n^0)$ plus
\begin{equation}
  P^n \bigl\{\bigl| \rave_{1}^n {(g-h)(x_i)}\bigr|
  > \delta \bigr\} \le \frac{\Vert g-h\Vert^2 }{\delta^2}
       \le \Tfrac{\Ds\delta^3}{\Ds\delta^2} = \delta
\end{equation}
This implies
\begingroup \mathsurround0em\arraycolsep0em \begin{eqnarray}
\noalign{\vspace{-\abovedisplayskip}\vspace{\abovedisplayshortskip}}
  & \Ds   \Nlimsup_{n}
  P^n \bigl\{| \tau_{n,g}^*- \tau_{n,h}^* | > \varepsilon\bigr\} \le
  8 \hspace{6\ei}\delta \hspace{9\ei} \frac{\varphi(0)}{\Vert h\Vert } +
  \delta & \\ \noalign{\noindent hence \nopagebreak}
\label{e.t.gtoht*gt*h}& \Ds   \lim_{g\to h}\Nlimsup_{n}
  P^n \bigl\{| \tau_{n,g}^*- \tau_{n,h}^* | > \varepsilon\bigr\}=0 &
\end{eqnarray}\endgroup
Observe that
      $ \limsup_{n} P^n \bigl\{ |\tau_n- \tau_{n,h}^*| >
         2 \hspace{6\ei}\varepsilon\bigr\}  $
does not depend on~$g$, therefore, may be bounded by
\begin{equation} \mathsurround0em\arraycolsep0em
\begin{array}{l} \Ds
   \lim_{g\to h} \Nlimsup_{n}
   P^n \bigl\{|\tau_n- \tau_{n,g}^*|>\varepsilon\bigr\} +
   P^n \bigl\{| \tau_{n,g}^*- \tau_{n,h}^* | > \varepsilon\bigr\}
  \hspace{-.75em} \\ \rule{0pt}{3.3ex} \Ds \hspace{2.25em}{}\le
   \lim_{g\to h} \Nlimsup_{n}
   P^n \bigl\{|\tau_n- \tau_{n,g}^*|>\varepsilon\bigr\} \\
\rule{0pt}{3.3ex} \Ds \hspace{4.5em}{}+
   \lim_{g\to h}\Nlimsup_{n}
   P^n \bigl\{| \tau_{n,g}^*- \tau_{n,h}^* | > \varepsilon\bigr\}
\end{array}\end{equation}
The upper bound equals zero by~(\ref{e.t.gtohPAn0})
and~(\ref{e.t.gtoht*gt*h}); thus,
\begin{equation} \label{e.t.tnth}
  \Nlimsup_{n} P^n \bigl\{ |\tau_n- \tau_{n,h}^*| >
         2 \hspace{6\ei}\varepsilon\bigr\}  = 0
\end{equation}
\par     It remains to prove that
\begin{equation} \label{e.t.tNPsum}
  \tau_{n,h}^*=\tau_n^\star + \Lo_{P^n}(n^0)
  \hspace{1.5em} \mbox{for} \hspace{.75em}
  \tau_n^\star = \hJc \bigl(\rave_{1}^n h(x_i) >
                \Vert h\Vert \hspace{6\ei} u_{\alpha}\bigr)
\end{equation}
But $ \rave_1^n h(x_i) $ is asy.~normal
  $ {\cal N} \bigl(0,\Vert h\Vert^2 \hspace{6\ei}\bigr )$ and
  $ {\cal N} \bigl( \hspace{6\ei}\Vert h\Vert^2,
               \Vert h\Vert^2 \hspace{6\ei}\bigr ) $
under~$P^n$, respectively~$P^n_{n,h}$, so that
\begin{equation}
  \Nlim_{n}\itg \tau_n^\star\,dP^n=\alpha\weg, \qquad
  \Nlim_{n}\itg \tau_n^\star\,dP^n_{n,h}=
  \Phi \bigl(-u_{\alpha}+ \Vert h\Vert\,\bigr)
\end{equation}
Thus, the uniqueness result of HR~(1994; Corollary~3.4.2\footnote{Note that
   $\sigma>0$ must be assumed in part~(b).}) applies to~$ (\tau_n^\star) $,
such that $ \tau_n^\star=\tau_{n,h}^* + \Lo_{P^n}(n^0) $. Altogether,
   (\ref{e.t.tnth}) and~(\ref{e.t.tNPsum}) imply~(\ref{e.t.1exp}).
\par     The converse, that~(\ref{e.t.uopt}) entails optimality,
is obvious, as all sequences~$(Q_n^n)$ in~$ H\cup K$ are contiguous
to~$(P^n)$.\qed \end{Bew}
\subsection{Invariant Tangent Cones and Spaces}
\paragraph{Rank Functionals}
For the symmetry problem on the real line,
    one-sample rank functionals~$R_{\varrho}$ are given by
\begin{equation}  \label{e.t.RFb}
     R_{\varrho}(Q) =
     2\itg_{0}^{\infty} \varrho\bigl( Q(x)-Q(-x)\bigr)\,Q(dx)
                    - \itg_{0}^1 \varrho \,d \lambda_0
\end{equation}
where  $Q(x)=Q \bigl((-\infty,x]\hspace{6\ei}\bigr)$,
       $\lambda_0$ denotes Lebesgue measure on~$(0,1)$,
and~$\varrho$ is some (scores) function in~$L_1(\lambda_0)$.
Then~$R_{\varrho}(Q)$ is defined for every~$Q\in {\cal M}_c$,
the set of all probabilities with continuous distribution functions.
Let~${\cal M}_{cs}$ denote the subset of all symmetric
  $ P\in {\cal M}_c $ (that is, $ P(-x)=1-P(x) \enskip \forall\,x>0 $).
Then $R_{\varrho}(P)=0$ for all $P\in {\cal M}_{cs}$.
A certain kind of asymmetry is defined through nonzero values
of the functional. If~$\varrho$ is nonnegative increasing,
then $R_{\varrho}(Q)\ge0$ for all positively asymmetric
  $ Q\in {\cal M}_{c} $ (that is, $ Q(-x)\le 1-Q(x) \enskip \forall\,x>0 $);
more generally,
  $ R_{\varrho}(Q'')\ge R_{\varrho}(Q') $ if $Q',Q''\in {\cal M}_{c}$,
  $ Q''(x)\le Q'(x)\enskip \forall\,x\in\R$.
\paragraph{Signed Linear Rank Statistics}
Linear rank statistics~$R_n$ are of the form
\begin{equation} \label{e.t.RSbn}
  R_n = \ave_{i=1}^n \zi(x_i)\,\varrho_{n}(r_{n,i}^{+})
\end{equation}
where $r_{n,i}^{+}$ denote the absolute ranks (rank~$|x_i|$ among
      $|x_1|,\ldots,|x_n|$), and~$\varrho_{n}(i)$ are some numbers (scores).
The weak condition used by H\'ajek and Sid\'ak~(1967; V.1.7) to
prove asy.~normality of~$R_n$ under~$P^n$ (in fact, asy.~linearity at~$P$) is
\begin{equation} \label{e.t.bnitob}
    \varrho_{n}([1+ns]) \longrightarrow \varrho(s)
    \qquad \mbox{in}\enskip   L_2(\lambda_0) \hspace{-2.25em}
\end{equation}
Given any $\varrho\in L_2(\lambda_0)$,
this condition is satisfied by the array
   $ \varrho_{n}(i) = \Ew \varrho(u_{n(i)}) $
(based on the order statistics~$u_{n(i)}$ of an i.i.d.~sample
   $ u_1,\ldots,u_n\sim \lambda_0 $), by the array
   $ \varrho_{n}(i) = n \itg_{I_n}\varrho\,d \lambda_0 $ with
   $ I_n = (\frac{i-1}{n},\frac{i}{n}) $, and the array
   $ \varrho_{n}(i) = \varrho(\Tfrac{i}{n+1}) $
(under a mild extra condition on~$\varrho$).
Then, for every $P\in {\cal M}_{cs}$, the sequence of rank
statistics~$ (R_n) $ is asy.~linear at~$P$ with influence
curve~$\kappa_P$,
\begingroup \mathsurround0em\arraycolsep0em
\begin{eqnarray} \label{e.t.asNRn}
     R_n &{}={}& \ave_{i=1}^n \kappa_P(x_i) + \Lo_{P^n}(1/\!\sqrt{n}\,) \\
\noalign{\noindent where\nopagebreak \vspace{-1ex}}\label{e.t.ICrank}
    && \kappa_P(x) = \zi(x)\, \varrho
    \bigl(\hspace{6\ei} 2 \hspace{9\ei}P(|x|)-1 \hspace{3\ei}\bigr)
\end{eqnarray}\endgroup
An alternative approach imposes bounds on the growth of the derivative(s)
of the scores function~$\varrho$; confer H\'ajek and Sid\'ak~(1967; VI.5.1).
These Chernoff--Savage conditions have successively been weakened and ensure
the asy.~normality of $ \sqrt{n}\,\bigl(R_n-R_{\varrho}(Q_n)\bigr)$,
even under noncontiguous alternatives~$(Q_n^n)$, with~$R_{\varrho}$ as
centering functional. Combining both sets of conditions,
differentiablity of~$R_{\varrho}$ at~$P\in {\cal M}_{cs}$ may be proved
as in HR~(1981\,a; Proposition~4.1).
Thus, at every $P\in {\cal M}_{cs}$, the functional~$R_{\varrho}$
is differentiable in the sense of~(\ref{e:i:Tdiff}) with influence
curve the same~$\kappa_P$ given by~(\ref{e.t.ICrank}).
\paragraph{Invariant Tangent Sets and Hypotheses}
Rank statistics~$R_n$ are not only distribution free under the null
hypothesis~${\cal M}_{cs}$ but also under suitably defined alternatives.
Let a family of sets~${\cal G}_P$, one for each $P\in {\cal M}_{cs}$,
be generated by some set ${\cal G}_0 \subset L_2(\lambda_0)$ such that
\begin{equation}\label{e.t.GP}
   {\cal G}_P = \bigl\{\, g_{P,q}\bigm|
                q\in {\cal G}_0 \hspace{6\ei}\bigr\}\weg, \qquad
   g_{P,q}(x) = \zi(x)\, q \bigl(\hspace{3\ei} 2 \hspace{6\ei}
                P(|x|)-1 \hspace{3\ei}\bigr)
\end{equation}
These sets~${\cal G}_P$, which obviously consist of odd functions,
are invariant in the sense that the composition
  $ {\cal G}_P \circ P^{-1} = \{\,g\circ P^{-1}\mid g\in {\cal G}_P\}$
with the pseudo-inverse $P^{-1}(s)=\inf \{x\in\R\mid P(x)\ge s\}$
is the same 
for all $P\in {\cal M}_{cs}$,
\begin{equation}
   {\cal G}_P \circ P^{-1}  = \bigl\{\, g_{0,q}\bigm|
                q\in {\cal G}_0 \hspace{6\ei}\bigr\}\weg, \qquad
   g_{0,q}(s) = \zi(s- \Tfrac{1}{2})\, q \bigl(\hspace{3\ei}|\hspace{3\ei}
   2 \hspace{6\ei}s-1 \hspace{3\ei}| \hspace{3\ei}\bigr)
\end{equation}
As $\itg g_{P,q} \,dP=0$ and $ \itg g^2_{P,q}\,dP= \itg q^2 \,d \lambda_0$,
the sets~${\cal G}_P$ may actually serve as tangent sets
   at~$P\in {\cal M}_{cs}$.
Moreover, the properties of~${\cal G}_0$ to be closed, convex, a cone,
a linear subspace of~$L_2(\lambda_0)$, respectively, are each inherited
to the sets~${\cal G}_P$ in~$L_2(P)$ for every $ P\in {\cal M}_{cs} $.
\begin{Rem}\rm Conversely, given any set~${\cal G}_{P_0}$ of odd tangents
at some $P_0\in {\cal M}_{cs}$, define
\begin{equation} \label{e.t.G0zuGP}
  {\cal G}_0 =
  \bigl\{\, q_{g}\bigm| g\in {\cal G}_{P_0}\hspace{6\ei}\bigr\}\weg, \qquad
  q_{g}(s)= g \bigl(P_0^{-1}(\Tfrac{1+s}{2})\bigr)
\end{equation}
Then this set~${\cal G}_0$, via~(\ref{e.t.GP}),
reproduces the given tangent set~${\cal G}_{P_0}$ at~$P_0$
and generates the following tangent sets~${\cal G}_P$
at other measures $P\in {\cal M}_{cs}$,
\begin{equation}
  {\cal G}_P= \bigl\{\, g \circ P_0^{-1}\circ P\bigm|
  g\in {\cal G}_{P_0}\,\bigr\}
\end{equation}
where $ g \circ P_0^{-1}\bigl(P(x)\bigr)=\zi(x)\hspace{15\ei}
           g \circ P_0^{-1}\bigl(P(|x|)\bigr)$~a.e.$P(dx)$.
Note that~$P_0^{-1}\circ P$ is odd and strictly increasing a.e.$P$.
For such tranformations applied to each~$x_i$, the vector of signs and
absolute ranks is (maximal) invariant.
\par
Positive shifts, for example, of some $P_0\in {\cal M}_{cs}$ which
has finite Fisher information of location and a Lebesgue density~$p_0$,
lead to the tangent cone generated by the function
  $ -(\hspace{6\ei}{\dot p}_0/p_0)=g_{P_0,q_0}$,
where $ q_0(s)= -(\hspace{6\ei}{\dot p}_0/p_0)\circ
                 P_0^{-1}(\frac{1+s}{2}) $,
and then $ g_{P_0,q_0} \circ P_0^{-1}\bigl(P(x)\bigr) =
           - \zi(x)\hspace{15\ei} (\hspace{6\ei}{\dot p}_0/p_0)\circ
           P_0^{-1} \bigl(P(|x|)\bigr) $~a.e.$P(dx)$.       \qed\end{Rem}
Now suppose that~${\cal G}_0$ is~(\ref{i.t.cone}) a convex cone,
or~(\ref{i.t.lin}) a linear space, in~$L_2(\lambda_0)$. \linebreak
For each $ P\in {\cal M}_{cs} $, let the hypotheses~$J^0_{P}$,
  $J_P$, $H_P$, and~$K_P$ about the rank functional~$R_{\varrho}$ over
the tangent set~${\cal G}_P$ be defined by~(\ref{e:t:J0})--(\ref{e:t:K}).
These hypotheses are invariant as they read
\begin{equation} \label{E:T:PinK/cone0}
  J^0_P: q=0\weg,\hspace{1em}
  J_P: {\langle \varrho|q\rangle}_{\hspace{-4\ei}0} = 0\weg,\hspace{1em}
  H_P: {\langle \varrho|q\rangle}_{\hspace{-4\ei}0}\le0\weg,\hspace{1em}
  K_P: t {\langle \varrho|q\rangle}_{\hspace{-4\ei}0} \ge c
\end{equation}
with reference to the tangent set~${\cal G}_P$ given by~(\ref{e.t.GP})
at~$P\in {\cal M}_{cs}$. In view of~(\ref{E:T:PinK/cone}),
representation~(\ref{E:T:PinK/cone0}) is a consequence of the following
equality of scalar products and norms in~$L_2(P)$ and~$L_2(\lambda_0)$,
respectively, for the tangents of form~(\ref{e.t.GP}),
\begin{equation} \label{e.t.normP=normL}
   {\langle \kappa_P|g\rangle}_{\hspace{-6\ei}P}
    = {\langle \varrho|q\rangle}_{\hspace{-4\ei}0}\weg,
\qquad {\Vert \kappa_P-g\Vert}^2_P= {\Vert \varrho -q\Vert }^2_0
\hspace{-1.5em}\end{equation}
\paragraph{Invariant Optimality of Rank Tests}
As another consequence of~(\ref{e.t.normP=normL}) we observe that
the approximation of~$\kappa_P$ by~$g\in {\cal G}_P$ is equivalent to
the approximation of~$\varrho$ by~$q\in {\cal G}_0$.
Therefore, 
the projection~$\hat{\kappa}_P$ of~$\kappa_P$ on~$\cl{\cal G}_P$
in~$L_2(P)$ is given in terms of the projection~$\hat{\varrho}$
of~$\varrho$ on~$\cl{\cal G}_0$ in~$L_2(\lambda_0)$,
\begin{equation}
    \hat{\kappa}_P(x)=\zi(x)\, \hat{\varrho}
    \bigl(\hspace{6\ei} 2 \hspace{9\ei}P(|x|)-1 \hspace{3\ei}\bigr)
\end{equation}
Then Theorem~\ref{t:t:1s-power} is in force and yields the optimal
asy.\ level~$\alpha$ test sequence~$(\hat{\tau}_{n,P})$
for $J^0_P$ vs.~$K_P$,
\begin{equation}
   \hat{\tau}_{n,P}=\hJc  \bigl(\rave_{i=1}^n
      \hat{\kappa}_P(x_i) > {\Vert \hat{\kappa}_P\Vert}_P
         \hspace{6\ei} u_{\alpha}\bigr)
\end{equation}
Now invoke any array of scores~$\hat{\varrho}_{n}(i)$ that,
via~(\ref{e.t.bnitob}), are connected to~$\hat{\varrho}$.
Employ the corresponding rank statistics~$\hat{R}_n$ to define
the rank tests
\begin{equation} \label{e.t.Rntestbhat}
   \hat{\tau}_{n}=\hJc  \bigl(\sqrt{n}\,\hat{R}_n >
   {\Vert \hat{\varrho}\Vert}_0 \hspace{6\ei} u_{\alpha}\bigr)
\end{equation}
independently of~$P\in {\cal M}_{cs}$.
Then, by~(\ref{e.t.asNRn}),~(\ref{e.t.ICrank})
for~$\hat{R}_n$ and~$\hat{\kappa}_P$,~$\hat{\varrho}$,
and by asy.~normality,
\begin{equation}
    \hat{\tau}_{n}=\hat{\tau}_{n,P}+ \Lo_{P^n}(n^0) \hspace{-.75em}
\end{equation}
for every $P\in {\cal M}_{cs}$.
Thus, the sequence~$(\hat{\tau}_n)$ of rank tests~(\ref{e.t.Rntestbhat})
is optimal for~$J^0_P$---if ${\cal G}_0=\lin{\cal G}_0$ even
for~$H_P$---against~$K_P$, according to Theorem~\ref{t:t:1s-power}.
\par     This optimality, in the two cases~(\ref{i.t.cone})
         ${\cal G}_0$ a convex cone, (\ref{i.t.lin})~${\cal G}_0$
a linear space, holds true for every $P\in {\cal M}_{cs}$.
\section{Confidence Limits}     \setcounter{equation}{0}\label{s.E}
Let~$P$ be any element of~${\cal P}$, with tangent set
    ${\cal G}\subset L_2(P)\cap \{1\}^{\perp}$, and some
    constant $c\in (0,\infty)$.
Similarly to the testing whether $ T(Q)\ge T(P)+c/\!\sqrt{n}\,$,
we now consider lower confidence limits $S_n-c/\!\sqrt{n}\,$ for
the value~$T(P)$; for example, the minimum amount of cash to be
kept on a business account.
Here and subsequently, the estimator sequence~$(S_n)$ may be any
sequence of estimates~$S_n$ at sample size~$n$.
It is desirable that~$S_n$ underestimate $T(P)+ c/\!\sqrt{n}\,$
with highest possible probability, under the i.i.d.~observations
  $x_1,\ldots,x_n\sim P$. This aim, however, is not well-defined,
as shown by arbitrary estimates~$S_n \le T(P)$.
Therefore, a side condition that also $S_n\ge T(P)$ with
sufficiently high probability must be imposed. In addition, to cut
out $S_n \equiv T(P)$, a local variation of~$P$ must be employed.
\subsection{Confidence Bounds For Lower and Upper Limits}
The following result requires some one-sided, respectively two-sided,
asymptotic median unbiasedness under the local perturbations~$P_{n,t,g}$
of~$P$ of kind~(\ref{e.t.Qpath}), and is of the type intended by
Pfanzagl and Wefelmeyer~(1982; Theorem~9.2.2).
\par
Qualitatively speaking, Theorem~\ref{t:e:conf.bd}(\ref{i.e.cone}) bounds
any `limit distribution function' of $\sqrt{n}\,\bigl(S_n-T(Q)\bigr)$
under~$Q=P$, subject to upper bound $1/2$ at the origin under
all~$Q=P_{n,t,\hat{\kappa}}$, on the positve half-line by that
of~${\cal N}\bigl(0, \Vert \hat{\kappa}\Vert^2 \hspace{6\ei}\bigr)$
from above.
In addition, Theorem~\ref{t:e:conf.bd}(\ref{i.e.lin}) bounds
such `limit distribution functions' under~$Q=P$, subject to the lower
bound~$1/2$ at the origin under all~$Q=P_{n,t,-\hat{\kappa}}$,
at the same time on the negative half-line by that
of~${\cal N}\bigl(0, \Vert \hat{\kappa}\Vert^2 \hspace{6\ei}\bigr)$
from below;
where $\hat{\kappa}=\tilde{\kappa},\bar{\kappa}$, respectively.
For best estimator accuracy, the limit distribution function should
be maximal on~$(0,\infty)$, and minimal on~$(-\infty,0)$.
In general, `limit distribution functions' need not exist nor need
they be normal.
\par
An estimator sequence~$(\hat S_n)$ is said to attain a confidence
upper bound
\begin{equation}
  \Nlimsup_{n} P^n \bigl\{ -t' <  \sqrt{n}\,\bigl(S_n-T(P)\bigr)
    < t'' \hspace{6\ei}\bigr\} \le \beta(t',t'')
\end{equation}
uniformly in $t'$,~$t''$, and with $\limsup_n$ replaced by~$\liminf_n$,
if~$(\hat S_n)$ itself satisfies the side conditions on the estimator
sequences~$(S_n)$ and in fact achieves 
\begin{equation} \label{e:e:conf.bd.ufo.achv.}
   \Nliminf_{n}\Ninf_{t',t''} \Bigl (
     P^n \bigl\{ -t' <  \sqrt{n}\,\bigl(\hat S_n-T(P)\bigr)
    < t'' \hspace{6\ei}\bigr\}- \beta(t',t'') \Bigr ) \ge0
\end{equation}
As for asy.~linear estimators, the reader is referred to the
beginning of Subsection~\ref{ss.e.reg}, where this kind of estimators
are introduced in more generality.
\begin{Thm}\sl \label{t:e:conf.bd}
Let~$(S_n)$ be any estimator sequence.
\begin{ABC} \item    \label{i.e.cone}
Suppose\/ {\rm ${\cal G}=\tilde{{\cal G}}$,} a convex cone.
Assume there exists some sequence of tangents $g_m\in \tilde{{\cal G}}$
such that $g_m\to \tilde{\kappa}$ in~$L_2(P)$ and, for every
convergent sequence $t_n\to t$ in\/~{\rm $(0,\infty)$,}
\begin{equation} \label{e.e.k.oomedS.ge0/cone}
  \Nliminf_{m}   \Nliminf_{n} P_{n,t_n,g_m}^n \bigl\{
   S_n \ge T(P_{n,t_n,g_m})  \bigr\} \ge \Tfrac{1}{2}
\end{equation}
Then, for every $t\in (0,\infty)$
      and every convergent sequence $t_n\to t$ in\/~{\rm $(0,\infty)$,}
\begin{equation} \label{e.e.low.conf.bd/cone}
  \Nlimsup_{n} P^n \bigl\{
  \sqrt{n}\,\bigl(S_n-T(P)\bigr) < t_n \bigr\} \le
  \Phi \Bigl(\frac{t}{\Vert \tilde{\kappa}\Vert }\Bigr)
\end{equation}
The upper bound\/~{\rm (\ref{e.e.low.conf.bd/cone})} is attained by
the asy.~linear estimator\/~{\rm $(\tilde{S}_n)$,}
\begin{equation}\label{e.e.tildeS}
  \sqrt{n}\,\bigl(\tilde{S}_n-T(P)\bigr)
  = \rave_{i=1}^n \tilde{\kappa}(x_i)   + \Lo_{P^n}(n^0)
\end{equation}
which achieves\/~{\rm (\ref{e.e.low.conf.bd/cone})}
with\/ {\rm $t_n=t$,} uniformly in\/~{\rm $-\infty\le t \le \infty$,}
and with $\limsup_n$ replaced by~$\liminf_n$.
\item  \label{i.e.lin}
Suppose\/ {\rm ${\cal G}=\bar{{\cal G}}$,} a linear space.
Assume there exist two sequences of tangents $g_m',g_m''\in \bar{{\cal G}}$
such that\/ {\rm $g_m'\to\bar{\kappa}$, $g_m''\to-\bar{\kappa}$} in~$L_2(P)$
and, for every convergent sequence $t_n\to t$ in\/~{\rm $(0,\infty)$,}
\begingroup \mathsurround0em\arraycolsep0em \begin{eqnarray}
\label{e.e.k.oomedS.ge0/lin}
  \Nliminf_{m}   \Nliminf_{n} P_{n,t_n,g_m'}^n \bigl\{
   S_n \ge T(P_{n,t_n,g_m'})  \bigr\} &{}\ge{}& \Tfrac{1}{2} \\
\label{e.e.k.oomedS.le0/lin}
  \Nliminf_{m}   \Nliminf_{n} P_{n,t_n,g_m''}^n \bigl\{
  S_n\le T(P_{n,t_n,g_m''}) \bigr\} &{}\ge{}& \Tfrac{1}{2}
\end{eqnarray}\endgroup
Then, for every $t',t''\in (0,\infty)$ and all sequences\/
      {\rm $t_n'\to t'$, $t_n''\to t''$} in\/~{\rm $(0,\infty)$,}
\begin{equation} \label{e.e.2sconf.bd/lin}
  \Nlimsup_{n} P^n \bigl\{
  -t_n' <   \sqrt{n}\,\bigl(S_n-T(P)\bigr) < t_n'' \hspace{6\ei}\bigr\}
  \le   \Phi \Bigl(\frac{t''}{\Vert \bar{\kappa}\Vert }\Bigr)   -
        \Phi \Bigl(-\frac{t'}{\Vert \bar{\kappa}\Vert }\Bigr)
\end{equation}
The upper bound\/~{\rm (\ref{e.e.2sconf.bd/lin})} is attained by the
asy.~linear estimator\/~{\rm $(\bar{S}_n)$,}
\begin{equation}\label{e.e.barS}
  \sqrt{n}\,\bigl(\bar{S}_n-T(P)\bigr) =
  \rave_{i=1}^n \bar{\kappa}(x_i) + \Lo_{P^n}(n^0)
\end{equation}
which achieves\/~{\rm (\ref{e.e.2sconf.bd/lin})}
with\/ {\rm $t_n'=t'$, $t_n''=t''$,} uniformly in\/
       {\rm $-\infty\le -t'<t''\le \infty$,} and with
            $\limsup_n$ replaced by~$\liminf_n$.
\end{ABC}\end{Thm}
\begin{Rem}\rm  [\,asymptotic median nonnegative, nonpositive\,]
\par     Conditions~(\ref{e.e.k.oomedS.ge0/cone}),
         (\ref{e.e.k.oomedS.ge0/lin}), and~(\ref{e.e.k.oomedS.le0/lin}),
respectively, mean that---in the iterated limit---the median
of~$\sqrt{n}\,\bigl(S_n-T(P_{n,t_n,g})\bigr)$ under~$ P_{n,t_n,g}^n $
for~$n$ large, $g\approx\tilde{\kappa}$, and $g\approx\bar{\kappa}$,
respectively, becomes~$\ge0$, and~$\le0$ for $g \approx -\bar{\kappa}$.
\par     Of course, if $\tilde{\kappa}\in \tilde{{\cal G}}$,
         respectively $\bar{\kappa}\in \bar{{\cal G}}$,
conditions~(\ref{e.e.k.oomedS.ge0/cone}) and
           (\ref{e.e.k.oomedS.ge0/lin}),~(\ref{e.e.k.oomedS.le0/lin})
are needed only for $g_m=\tilde{\kappa}$, respectively
    for $ g_m'=\bar{\kappa} $ and~$g_m''=-\bar{\kappa}$.
\par    Conditions~(\ref{e.e.k.oomedS.ge0/cone}),
        (\ref{e.e.k.oomedS.ge0/lin}) and~(\ref{e.e.k.oomedS.le0/lin}),
respectively, are ensured by asymptotic median nonnegativity and
nonpositivity, respectively, for every fixed tangent in the corresponding
tangent set~${\cal G}$, in the sense of~(\ref{e.e.oomedS.ge0.g/all})
        and~(\ref{e.e.oomedS.le0.g/all}) below.\qed         \end{Rem}
\begin{Bew}\enskip
We start the derivation of the bounds simultaneously in both cases:
\par     Fix any $g\in {\cal G}$ such that $ \langle \kappa|g \rangle \ne0$,
      any sequence $t_n\to t$ in~$(0,\infty)$, and put $P_n= P_{n,t_n,g}$.
Expansion~(\ref{e.t.expT/twn}), by~(\ref{e:i:Tdiff}), holds uniformly
on $t$-compacts, so
\begin{equation} \label{e.e.sntkg}
    \sqrt{n}\,\bigl(\hspace{6\ei}T(P_{n})-T(P)\bigr) =  t \hspace{6\ei}
    \langle \kappa|g\rangle + \Lo(n^0) = s_n \langle \kappa|g\rangle
\end{equation}
for some suitable other sequence $s_n=s_{n,t_n,g}\to t$. Thus, we obtain
\begin{equation} \label{e.e.STR}
   \sqrt{n}\,\bigl( S_{n}-T(P_n)\bigr)= R_n-s_n \langle \kappa|g\rangle
   \hspace{1.5em}\mbox{for}\hspace{1.5em}
   R_n=\sqrt{n}\,\bigl( S_{n}-T(P)\bigr)
\end{equation}
Also the loglikelihood expansion~(\ref{E:T:llh}) for fixed~$g$,
due to~(\ref{e:i:pathdiff}), holds uniformly on~$t$-compacts.
Therefore, and by mutual contiguity of~$(P_n^n)$ and~$(P^n)$,
\begin{equation}\label{E:e:llh.tn/+/-}
  \log \frac{dP^n}{dP_{n}^n} =
  - \log \frac{dP_{n}^n}{dP^n} + \Lo_n' =
    - t \rave_{i=1}^n g(x_i)
  + \Tfrac{1}{2}  t^2 \hspace{6\ei}{\Vert g\Vert}^2 + \Lo_n''
\end{equation}
where $ \Lo_n'$,~$ \Lo_n''$ each are some $\Lo_{P^n}(n^0)$.
By HR~(1994; Proposition~2.2.12 and Corollary~3.4.2\,a),
the asy.~power of any test sequence~$(\tau_n)$ under~$(P^n)$,
subject to asy.~level~$\alpha$ under~$(P_n^n)$, is bounded
by~$\Phi \bigl(-u_{\alpha}+t \Vert g\Vert \,\bigr)$.
\par     Applying this bound to the sequence of tests 
\begingroup \mathsurround0em\arraycolsep0em \begin{eqnarray}
\label{e.e.taunRn}   \tau_n &{}={}&
  \Jc ( R_n< s_n \langle \kappa|g\rangle \hspace{6\ei}) = \tau_{n,t_n,g} \\
\noalign{\noindent and their asy.~level\nopagebreak}
\label{e.e.level.ag}   \alpha_g &{}={}& \Nlimsup_{n} P_n^n \bigl\{
            R_n< s_n \langle \kappa|g\rangle \bigr\} \hspace{-1.5em}\\
\noalign{\noindent we obtain\nopagebreak} \label{e.e.1s.pre-bd}
   \Phi \bigl(-u_{\alpha_g}+t \Vert g\Vert \hspace{12\ei}\bigr) &{}\ge{}&
   \Nlimsup_{n} P^n \bigl\{ R_n< s_n \langle \kappa|g\rangle
   \bigr\}     \hspace{-1.5em}
\end{eqnarray}\endgroup
\begin{ABC}\item[{\rfi{\ref{i.e.cone}}}] 
Observe that, by condition~(\ref{e.e.k.oomedS.ge0/cone}),
  as $g=g_m\in \tilde{{\cal G}}$ tends to~$\tilde{\kappa}$ in~$L_2(P)$,
\begin{equation}
   \limsup\alpha_{g} \le \Tfrac{1}{2}\weg, \quad
   \mbox{hence} \quad \liminf u_{\alpha_{g}} \ge0
\end{equation}
Therefore, given $\delta\in (0,1)$, one can choose
  $g=g_m\in \tilde{{\cal G}}$ so close to~$\tilde{\kappa}$ that
\begin{equation} \label{e.e.uagSkg}
   -u_{\alpha_g}+t \hspace{6\ei}\Vert g\Vert \le
    t \hspace{6\ei}\Vert \tilde{\kappa}\Vert + \delta
   \qquad \mbox{and}\qquad    s_n \langle \kappa|g\rangle
   \ge (t_n- \delta)\Vert \tilde{\kappa}\Vert^2
\end{equation} 
eventually. Then~(\ref{e.e.1s.pre-bd}) implies that
\begingroup \mathsurround0em\arraycolsep0em \begin{eqnarray}
  \Nlimsup_{n} P^n \bigl\{  R_n< (t_n- \delta)
  \Vert \tilde{\kappa}\Vert^2  \bigr\} &{}\le{}&   \Phi \bigl(
  \hspace{6\ei}  t \hspace{6\ei}\Vert \tilde{\kappa}\Vert +
  \delta\hspace{4\ei}\bigr) \\
\noalign{\noindent hence\nopagebreak} \label{e.e.1s.bd-d}
  \Nlimsup_{n} P^n \bigl\{  R_n< t_n
  \Vert \tilde{\kappa}\Vert^2  \bigr\} &{}\le{}&  \Phi \bigl( \hspace{6\ei}
  t \hspace{6\ei}\Vert \tilde{\kappa}\Vert + \delta \hspace{6\ei}
  \Vert \tilde{\kappa}\Vert + \delta\hspace{4\ei}\bigr)
\end{eqnarray}\endgroup
where assumption~(\ref{e.e.k.oomedS.ge0/cone}) has been used
for the shifted sequence~$t_n+ \delta$. Once more
using~(\ref{e.e.k.oomedS.ge0/cone}) for the rescaled
sequence~$t_n \big/\Vert \tilde{\kappa}\Vert^2 $,
bound~(\ref{e.e.low.conf.bd/cone}) follows from~(\ref{e.e.1s.bd-d}),
if we let $\delta\to0$.
\item[{\rfi{\ref{i.e.lin}}}] 
Starting from assumption~(\ref{e.e.k.oomedS.ge0/lin}),
the proof~(\ref{i.e.cone}) establishes the bound 
\begin{equation} \label{e.e.2s.confbd.1}
  \limsup_{n\to \infty} P^n \bigl\{  R_n< t_n'' \hspace{6\ei}
  \Vert \bar{\kappa}\Vert^2  \bigr\} \le \Phi \bigl( \hspace{6\ei} t''
  \hspace{6\ei}\Vert \bar{\kappa}\Vert \hspace{12\ei}\bigr)
\end{equation} 
for every $t''\in (0,\infty)$ and every convergent sequence
          $t_n''\to t''$ in~$(0,\infty)$.
\par     In addition, given $g\in \bar{{\cal G}}$ and another sequence
       $t_n'\to t'$ in~$(0,\infty)$, abbreviate~$P_{n,t_n',g}$ by~$Q_n$
    and choose~$r_n=r_{n,t_n',g}\to t'$ to satisfy~(\ref{e.e.sntkg})
    for~$(t_n')$.
\par     Then, like~(\ref{e.e.1s.pre-bd}) has been obtained for the
tests~(\ref{e.e.taunRn}), we conclude that
\begingroup \mathsurround0em\arraycolsep0em \begin{eqnarray}
\label{e.e.2s.pre-bd}
   \Phi \bigl(-u_{\beta_g}+t' \hspace{6\ei}
   \Vert g\Vert \hspace{12\ei}\bigr)    &{}\ge{}&
   \Nlimsup_{n} P^n \bigl\{
   R_n> r_n \langle \kappa|g\rangle \bigr\}     \hspace{-2.25em}\\
\noalign{\noindent using the tests\nopagebreak} \label{e.c.ups.tests}
  \upsilon_n &{}={}& \Jc ( R_n>r_n \langle \kappa|g\rangle \hspace{6\ei}) \\
\noalign{\noindent and their asy.~level\nopagebreak}
  \beta_g &{}={}& \Nlimsup_{n} Q_n^n \bigl\{
           R_n>r_n \langle \kappa|g\rangle \bigr\}  \hspace{-2.25em}
\end{eqnarray}\endgroup
By condition~(\ref{e.e.k.oomedS.le0/lin}),
as $g=g_m''\in \bar{{\cal G}}$ tends to~$-\bar{\kappa}$ in~$L_2(P)$,
\begin{equation}
     \limsup \beta_g \le \Tfrac{1}{2}\weg, \quad \mbox{hence} \quad
     \liminf u_{\beta_g} \ge0
\end{equation}
Therefore, given $\delta\in (0,1)$, we may choose
  $g=g_m''\in \bar{{\cal G}}$ so close to~$-\bar{\kappa}$ that
\begin{equation}
   -u_{\beta_g}+t' \hspace{6\ei}\Vert g\Vert \le
   t' \hspace{6\ei} \Vert \bar{\kappa}\Vert + \delta
   \qquad \mbox{and}\qquad    r_n \langle \kappa|g\rangle
   \le -(t_n'- \delta)\Vert \bar{\kappa}\Vert^2
\end{equation} 
eventually. Then~(\ref{e.e.2s.pre-bd}) implies that,
for each $\delta\in (0,1)$,
\begingroup \mathsurround0em\arraycolsep0em \begin{eqnarray}
  \Nlimsup_{n} P^n \bigl\{  R_n > - (t_n'- \delta)
  \Vert \bar{\kappa}\Vert^2  \bigr\} &{}\le{}& \Phi \bigl( \hspace{6\ei}
   t' \hspace{6\ei}\Vert \bar{\kappa}\Vert + \delta\hspace{4\ei}\bigr) \\
\noalign{\noindent hence\nopagebreak} \label{e.e.2s.confbd.20}
  \Nlimsup_{n} P^n \bigl\{  R_n > -t_n' \hspace{6\ei}
  \Vert \bar{\kappa}\Vert^2  \bigr\} &{}\le{}&  \Phi \bigl( \hspace{6\ei}
   t' \hspace{6\ei}\Vert \bar{\kappa}\Vert \hspace{12\ei}\bigr) \\
\noalign{\noindent that is,\nopagebreak} \label{e.e.2s.confbd.2}
 \Nliminf_{n} P^n \bigl\{  R_n \le -t_n' \hspace{6\ei}
  \Vert \bar{\kappa}\Vert^2  \bigr\} &{}\ge{}&  \Phi \bigl( -
   t' \hspace{6\ei}\Vert \bar{\kappa}\Vert \hspace{12\ei}\bigr)
\end{eqnarray}\endgroup
As
\begin{displaymath}\mathsurround0em\arraycolsep0em \begin{array}{l} \Ds
  \Nlimsup_{n} P^n \bigl\{ -t_n' \hspace{6\ei}
  \Vert \bar{\kappa}\Vert^2     < R_n < t_n'' \hspace{6\ei}
  \Vert \bar{\kappa}\Vert^2\bigr\} \le {}\\
\Ds \rule{0pt}{3ex}   \hspace{2.25em}   \Nlimsup_{n} P^n
  \bigl\{  R_n< t_n'' \hspace{6\ei}\Vert \bar{\kappa}\Vert^2  \bigr\}
  - \Nliminf_{n} P^n \bigl\{  R_n \le -t_n'  \hspace{6\ei}
     \Vert \bar{\kappa}\Vert^2  \bigr\} \hspace{-.75em}
\end{array} \end{displaymath}
bound~(\ref{e.e.2sconf.bd/lin}) follows from~(\ref{e.e.2s.confbd.1})
and~(\ref{e.e.2s.confbd.2}).                                    \end{ABC}
\vspace{\parskip}\par
We shall check attainment of the bounds simultaneously in both cases:
\par     \noindent
The asy.~linearity~(\ref{e.e.tildeS}) and~(\ref{e.e.barS}) entail
asy.~normality under~$P^n$,
\begin{equation} \label{e.e.asyNhatS}
   \Bigl \lgroup \sqrt{n}\,\bigl(\hat{S}_n-T(P)\bigr) \Bigr \rgroup
   \bigl( P^n\bigr) \gwto
   {\cal N}\bigl( 0, \Vert \hat{\kappa}\Vert^2 \hspace{6\ei}\bigr)
\end{equation}
for $\hat{S}_n=\tilde{S}_n$ with $\hat{\kappa}=\tilde{\kappa}$,
respectively for $\hat{S}_n=\bar{S}_n$ with $\hat{\kappa}=\bar{\kappa}$.
It follows that
\begin{equation} \label{e.e.2s.conf.bd.at/both}
   \Nlim_{n} P^n \bigl\{  -t' <
   \sqrt{n}\,\bigl(\hat{S}_n-T(P) \bigr) < t'' \hspace{6\ei}\bigr\} =
   \Phi \Bigl(\frac{t''}{\Vert \hat{\kappa}\Vert }\Bigr)
   -  \Phi \Bigl(-\frac{t'}{\Vert \hat{\kappa}\Vert }\Bigr)
\end{equation}
uniformly in $-\infty\le -t'<t''\le \infty$, in both cases.
\par       Verification of the regularity
           condition~(\ref{e.e.k.oomedS.ge0/cone}) for~$(\tilde{S}_n)$,
and of conditions (\ref{e.e.k.oomedS.ge0/lin})
              and~(\ref{e.e.k.oomedS.le0/lin}) for~$(\bar{S}_n)$,
is postponed to Subsection~\ref{sss.cr.reg}.
\qed \end{Bew}
\begin{Rem}\rm   \label{r:e:positivepart}
In Theorem~\ref{t:e:conf.bd}(\ref{i.e.cone}), the {\sl upper\/} bound
  $ \Phi \bigl(t /\Vert \tilde{\kappa}\Vert \hspace{6\ei}\bigr) $
on~$(0,\infty)$ given by~(\ref{e.e.low.conf.bd/cone}),
for the `limit distribution function'
of $\sqrt{n}\,\bigl(S_n-T(Q)\bigr)$ under~$Q=P$, does not extend
to a {\sl lower\/} bound on~$(-\infty,0)$, whereas the bound
  $ \Phi \bigl(t /\Vert \bar{\kappa}\Vert \hspace{6\ei}\bigr) $
does in Theorem~\ref{t:e:conf.bd}(\ref{i.e.lin}).
For example, given any $a\in (0,\infty)$,
consider the following modification~$(\breve{S}_n)$ of~$(\tilde{S}_n)$,
\begin{equation}
   \breve{S}_n= \tilde{S}_n\lor \bigl(\hspace{6\ei}
                 T(P)-a/\!\sqrt{n}\hspace{15\ei}\bigr)
\end{equation}
Then, if $g\in \tilde{{\cal G}}$ is such that
   $ \langle \kappa|g\rangle\ge0 $, and $t_n\to t$ in~$(0,\infty)$,
it holds that, eventually, $ T(P_{n,t_n,g}) \ge T(P)-a/\!\sqrt{n}\,$.
Using the asymptotic median nonnegativity~(\ref{e.e.oomedS.ge0.g/all})
of~$(\tilde{S}_n)$ to be proved in Subsection~\ref{sss.cr.reg},
we obtain that, eventually,
\begin{equation} 
      P_{n,t_n,g}^n \bigl\{ \breve{S}_n \ge T(P_{n,t_n,g}) \bigr\} =
      P_{n,t_n,g}^n \bigl\{ \tilde{S}_n \ge T(P_{n,t_n,g}) \bigr\}
      \ge \Tfrac{1}{2} + \Lo(n^0)
\end{equation}
Under~$P$, however, since
  $ \sqrt{n}\,\bigl(\breve{S}_n - T(P)\bigr) = (-a)\lor
    \sqrt{n}\,\bigl(\tilde{S}_n - T(P)\bigr) $, we have
\begin{equation} \mathsurround0em\arraycolsep0em
\hspace{-1.5em}\begin{array}{rl@{\hspace{2em}\mbox{if}\hspace{.5em}}l}
\Ds  P^n \bigl\{ \sqrt{n}\,\bigl(\breve{S}_n - T(P)\bigr)\le t \bigr\}
     = {} & \Ds 0 & \Ds t<-a\\
\rule{0pt}{3.3ex} \Ds \ggrto {} &
\Ds  \Phi \Bigl(\frac{t}{\Vert \tilde{\kappa}\Vert}\Bigr) & \Ds t\ge-a
\end{array}\end{equation}
The choice $a=0$ is possible if the asymptotic median~$\ge0$
condition~(\ref{e.e.oomedS.ge0.g/all}) is required,
instead of for $ \langle \kappa|g\rangle\ge0 $,
only for $ \langle \kappa|g\rangle>0 $,
which suffices for~(\ref{e.e.k.oomedS.ge0/cone}).\qed \end{Rem}
\subsection{Uniqueness of Efficient Estimators}            \label{ss.e.u}
In the setup of Theorem~\ref{t:e:conf.bd}(\ref{i.e.lin}), the optimal
estimates~$\bar{S}_n$ defined by~(\ref{e.e.barS}) are unique, up to terms
tending stochastically to zero under~$(P^n)$.
In the setup of Theorem~\ref{t:e:conf.bd}(\ref{i.e.cone}), on the contrary,
only the positive part $\bigl(\tilde{S}_n-T(P)\bigr)_+$ of the optimal
estimates~(\ref{e.e.tildeS}) centered at~$T(P)$ will be asymptotically
unique; confer Remark~\ref{r:e:positivepart} for an example.
\begin{Prop}\sl \label{p.e.u}
Let~$(\breve{S}_n)$ and~$(\hat{S}_n)$ be two estimator sequences.
\begin{ABC}\item[{\rfi{\ref{i.e.cone}}}]
In the case of a convex tangent cone\/~{\rm $\tilde{{\cal G}}$,}
suppose~$(\breve{S}_n)$ satisfies
  condition\/~{\rm (\ref{e.e.k.oomedS.ge0/cone})} and achieves
the confidence bound\/~{\rm (\ref{e.e.low.conf.bd/cone}),}
with $\limsup_n$ replaced by~$\liminf_n$. Then necessarily
\begin{equation} \label{e.e.u.Shat.cone}
\mathsurround0em\arraycolsep0em \begin{array}{r@{{}={}}l}
\Ds   \sqrt{n}\, \bigl(\breve{S}_n-T(P)\bigr)_+ + \bvLo_{P^n}(n^0)
    & \Ds \bigl(\rave_{i=1}^n \tilde{\kappa}(x_i)\bigr)_+  \\
\rule{0pt}{3ex} & \Ds
   \sqrt{n}\, \bigl(\tilde{S}_n-T(P)\bigr)_+ + \tlLo_{P^n}(n^0)
\end{array}\end{equation}
Conversely, form\/~{\rm (\ref{e.e.u.Shat.cone})} of~$(\breve{S}_n)$
implies\/~{\rm (\ref{e.e.oomedS.ge0.g/all})} and achievement of
bound\/~{\rm (\ref{e.e.low.conf.bd/cone}),} 
uniformly in\/~{\rm $-\infty\le t_n=t \le \infty$,}
and with $\limsup_n$ replaced by~$\liminf_n$.
\item[{\rfi{\ref{i.e.lin}}}]              In the case of a linear tangent
space\/~{\rm $\bar{{\cal G}}$,} assume~$(\hat{S}_n)$ satisfies
  conditions\/~{\rm (\ref{e.e.k.oomedS.ge0/lin})}
         and\/~{\rm (\ref{e.e.k.oomedS.le0/lin}),}
and achieves the confidence bound\/~{\rm (\ref{e.e.2sconf.bd/lin}),}
with $\limsup_n$ replaced by~$\liminf_n$. Then necessarily
\begin{equation} \label{e.e.u.Shat.lin}
\mathsurround0em\arraycolsep0em \begin{array}{r@{{}={}}l}
\Ds   \sqrt{n}\, \bigl(\hat{S}_n-T(P)\bigr) + \htLo_{P^n}(n^0)
   & \Ds \rave_{i=1}^n \bar{\kappa}(x_i)  \\
\rule{0pt}{3ex} & \Ds
   \sqrt{n}\, \bigl(\bar{S}_n-T(P)\bigr) + \brLo_{P^n}(n^0)
\end{array}\end{equation}
Conversely, if $(\hat{S}_n)$ is of form\/~{\rm (\ref{e.e.u.Shat.lin}),}
then it satisfies\/~{\rm (\ref{e.e.oomedS.ge0.g/all}),
                         (\ref{e.e.oomedS.le0.g/all}),} and
achieves bound\/~{\rm (\ref{e.e.2sconf.bd/lin}),}
with\/ {\rm $t_n'=t'$, $t_n''=t''$,} uniformly in\/
       {\rm $-\infty\le -t'<t''\le \infty$,} and with
            $\limsup_n$ replaced by~$\liminf_n$.
\end{ABC}\end{Prop}
\begin{Bew}\enskip
The proof draws on the proofs to Proposition~\ref{p.t.u}
and Theorem~\ref{t:e:conf.bd}:
\par     In case~(\ref{i.e.cone}), let~$(\breve{S}_n)$
satisfy~(\ref{e.e.k.oomedS.ge0/cone}) and achieve
bound~(\ref{e.e.low.conf.bd/cone}) such that,
for every constant sequence $t_n=t\in (0,\infty)$,
\begingroup \mathsurround0em\arraycolsep0em \begin{eqnarray}
& \Ds  \hspace{-1.5em}  \Nlimsup_{m}   \Nlimsup_{n} P_{n,t,g_m}^n
  \bigl\{   \breve{S}_n < T(P_{n,t,g_m})  \bigr\} \le \Tfrac{1}{2}
  \hspace{3.75em}\mbox{(\ref{e.e.k.oomedS.ge0/cone})}\hspace{-5em}
  \nonumber & \\
\noalign{\noindent and\nopagebreak} \label{e.c.u.att.cone}
& \Ds  \hspace{-1.5em} \Nliminf_{n} P^n \bigl\{ \breve{R}_n <
  t \hspace{6\ei} \Vert \tilde{\kappa}\Vert^2 \hspace{6\ei} \bigr\}
  \ge \Phi (t \hspace{6\ei}\Vert \tilde{\kappa}\Vert \hspace{6\ei})
\hspace{-.75em} & \end{eqnarray}\endgroup
where $ \breve{R}_n=\sqrt{n}\,\bigl(\breve{S}_n-T(P)\bigr)$.
Fix some $t_n=t\in (0,\infty)$ and any $\delta_a\in (0,1)$.
Choose $\delta\in (0,t)$ small enough and then
       $g=g_m\in \tilde{{\cal G}}$ so close to~$\tilde{\kappa}$ that
\begin{equation}
  \alpha_g < \Tfrac{1}{2}+\delta_a \weg, \qquad
  \Phi \bigl((t - \delta)\Vert \tilde{\kappa}\Vert \hspace{6\ei}\bigr)
  + \delta_a \ge \Phi(t \hspace{6\ei}\Vert g\Vert \hspace{6\ei})  \ge
  \Phi (t \hspace{6\ei}\Vert \tilde{\kappa}\Vert \hspace{6\ei}) - \delta_a
\end{equation}
and such that (\ref{e.e.uagSkg}) is fulfilled, too.
Recall (\ref{e.e.sntkg}), (\ref{e.e.STR}), and~(\ref{e.e.level.ag}).
Then
\begingroup \mathsurround0em\arraycolsep0em \begin{eqnarray}
& \Ds \hspace{-1.5em} \begin{array}{r@{{}\ge{}}l}
\Ds  \Nliminf_n P^n \bigl\{
     \breve{R}_n< s_n \langle \kappa|g\rangle \bigr\} &
\Ds  \Nliminf_n P^n \bigl\{
     \breve{R}_n< (t- \delta)\Vert \tilde{\kappa}\Vert^2
     \hspace{6\ei} \bigr\} \\ \rule{0pt}{3.5ex} & \Ds
   \Phi \bigl((t - \delta)\Vert \tilde{\kappa}\Vert \hspace{6\ei}\bigr)
   \ge \Phi(t \hspace{6\ei}\Vert g\Vert \hspace{6\ei}) - \delta_a
\end{array} & \\ \noalign{\noindent while\nopagebreak}
& \Ds \hspace{-1em}   \Nlimsup_{n} P_{n,t,g}^n \bigl\{
   \breve{R}_n< s_n \langle \kappa|g\rangle
   \bigr\}  =\alpha_g< \Tfrac{1}{2}+\delta_a &
\end{eqnarray}\endgroup
Therefore, the tests $ \tau'_{n,g} = 1- \tau_{n,t,g} = \Jc ( \breve{R}_n
                       \ge s_n \langle \kappa|g\rangle \hspace{6\ei}) $
given by~(\ref{e.e.taunRn}) satisfy
\begingroup \mathsurround0em\arraycolsep0em \begin{eqnarray}
& \Ds  \Nlimsup_n \itg \tau'_{n,g}\,dP^n \le
  \alpha'_g + \delta_a \le \alpha' + 2 \hspace{6\ei}\delta_a & \\
& \Ds  \Nliminf_n \itg \tau'_{n,g}\,dP_{n,t,g}^n \ge
  \Tfrac{1}{2} - \delta_a &
\end{eqnarray}\endgroup
where $ \alpha' =
         \Phi(-t \hspace{6\ei}\Vert \tilde{\kappa}\Vert \hspace{6\ei})$,
        $ \alpha'_g= \Phi(-t \hspace{6\ei}\Vert g\Vert \hspace{6\ei})$,
and so       $u_{\alpha'}=t \hspace{6\ei}\Vert \tilde{\kappa}\Vert $,
             $u_{\alpha'_g}=t \hspace{6\ei}\Vert g\Vert $.
Replacing~$\alpha$ and~$\beta_g$ in~(\ref{e.t.NPTg.abg}) by~$\alpha'_g$
and $\beta'_g = \Phi \bigl(-u_{\alpha'_g} + t \hspace{6\ei}
                \Vert g\Vert\hspace{6\ei}\bigr) = 1/2 = \beta' $,
respectively, (\ref{e.t.taun.anbndelta})~is satisfied by the
tests~$\tau'_{n,g}$ and leeway~$\delta_a$, in the place of~$\tau_n$
and~$\delta$ there. Via~(\ref{e.t.gtohPAn0}) and~(\ref{e.t.gtoht*gt*h}), we
reach~(\ref{e.t.tnth}). Taking already the asymptotic
equivalence~(\ref{e.t.tNPsum}) into account, where
     $ \Vert \tilde{\kappa}\Vert \hspace{6\ei} u_{\alpha'}
       = t \hspace{6\ei}\Vert \tilde{\kappa}\Vert^2 $,
and the fact that the present tests are all nonrandomized, we thus obtain
\begin{equation}\label{e.c.u.taungp/taunstar}
   \lim_{g\to \tilde{\kappa}} \Nlimsup_n
   P^n (\tau'_{n,g}\ne \tau_{n}^\star )  = 0 \weg,  \qquad
   \tau_{n}^\star = \hJc \bigl( \rave_{1}^n \tilde{\kappa}(x_i) \ge
   t \hspace{6\ei}\Vert \tilde{\kappa}\Vert^2 \hspace{6\ei}\bigr)
\end{equation}
The tests  $ \tau'_{n,g} = \Jc ( \breve{R}_n \ge s_n
              \langle \kappa|g\rangle \hspace{1\ei}) $
may be compared with
           $ \tau'_n = \Jc ( \breve{R}_n \ge t \hspace{6\ei}
             \Vert \tilde{\kappa}\Vert^2 \hspace{3\ei}) $.
In~(\ref{e.c.u.att.cone}),
    $ P^n \bigl( \breve{R}_n < t \hspace{6\ei}
      \Vert \tilde{\kappa}\Vert^2 \hspace{3\ei} \bigr) $
must actually converge to
    $ \Phi (t \hspace{6\ei}\Vert \tilde{\kappa}\Vert \hspace{3\ei}) $,
    and $s_n\to t$.
Therefore, employing the modulus~$\omega_{\Phi}$ of uniform continuity
of~$\Phi$, we obtain
\begin{equation}
  \Nlimsup_n   P^n (\tau'_{n,g}\ne\tau'_{n}\hspace{3\ei})
  \le \omega_{\Phi}(t \hspace{6\ei}|\lambda_g|\hspace{4\ei})\weg, \qquad
  \lambda_g= 
  \frac{\langle \kappa|g\rangle}{\Vert \tilde{\kappa}\Vert} 
             - \Vert \tilde{\kappa}\Vert
\end{equation}
As $\lim_{g\to \tilde{\kappa}}\lambda_g=0$, it follows that
\begin{equation} \label{e.c.u.taungp/tnp}
  \lim_{g\to \tilde{\kappa}}  \Nlimsup_n
  P^n (\tau'_{n,g}\ne \tau'_{n}) = 0
\end{equation}
Using the triangle inequality, we deduce from
(\ref{e.c.u.taungp/taunstar}) and~(\ref{e.c.u.taungp/tnp}) that
\begin{equation}
  \Nlim_n P^n \bigl\{\hspace{6\ei} \Jc ( \breve{R}_n \ge t \hspace{4\ei})
  \ne \Jc \bigl( \rave_{1}^n \tilde{\kappa}(x_i) \ge t \bigr) \bigr\} =0
\end{equation}
for every~$t\in \nolinebreak(0,\infty)$.
Because $\rave_{1}^n \tilde{\kappa}(x_i)$ is tight under~$(P^n)$,
the difference between the positive parts of~$\breve{R}_n$
and~$\rave_{1}^n \tilde{\kappa}(x_i)$ must converge to zero
in $P^n$-probability; confer HR~(1981\,b), fact (3.12)--(3.13).
Thus (\ref{e.e.u.Shat.cone})~is proved.
\par     In case~(\ref{i.e.lin}), we may now continue the same way as
part~(\ref{i.e.lin}) of the proof to Theorem~\ref{t:e:conf.bd} proceeds
after part~(\ref{i.e.cone}). From (\ref{e.t.NPTg.abg})
and~(\ref{e.t.taun.anbndelta}) onwards, plug the tests
  $ \upsilon'_{n,g} = 1- \upsilon_{n,t',g} = \Jc
    ( \hat{R}_n \le r_n \langle \kappa|g\rangle \hspace{1\ei}) $ given
by~(\ref{e.c.ups.tests}) in the proof of Proposition~\ref{p.t.u}.
Letting~$g$ tend to~$-\bar{\kappa}$, one similarly obtains that
\begin{equation}
  \Nlim_n P^n \bigl\{ \hspace{6\ei}\Jc ( \hat{R}_n \le -t \hspace{4\ei})
      \ne \Jc \bigl(  \rave_{1}^n \bar{\kappa}(x_i) \le -t \bigr) \bigr\} =0
\end{equation}
for every~$t\in \nolinebreak(0,\infty)$. This implies that
also the difference between the negative parts of~$\hat{R}_n$
and~$\rave_{1}^n \bar{\kappa}(x_i)$ must go to zero in $P^n$-probability,
hence~(\ref{e.e.u.Shat.lin}).
\par     As for the converse, which is obvious in case~(\ref{i.e.lin}),
observe in case~(\ref{i.e.cone}) that, for some stochastic term
         $\Lo_n=\Lo_{P^n}(n^0)$, and for every $t\in (0,\infty)$,
\begin{equation} \mathsurround0em\arraycolsep0em
\begin{array}{r@{{}={}}l}
\Ds  P^n (\breve{R}_n<t \hspace{4\ei}) & \Ds
     P^n (\breve{R}_n^+<t \hspace{4\ei}) =
     P^n (\tilde{R}_n^+ < t +\Lo_n) \\
\rule{0pt}{3ex} & \Ds    P^n (\tilde{R}_n^+ < t \hspace{4\ei}) + \Lo(n^0) =
                         P^n (\tilde{R}_n < t \hspace{4\ei}) + \Lo(n^0)
\end{array}\end{equation}
where the third equality is true because the
limit~$\Phi (t/ \Vert \tilde{\kappa}\Vert \hspace{6\ei})$
is continuous in~$t$. Thus, $(\breve{S}_n)$ inherits the optimality
from~$(\tilde{S}_n)$.
\par     Verification of the regularity conditions is postponed
to Subsection~\ref{sss.cr.reg}.
\qed\end{Bew}
\subsection{Regularity of Efficient Estimators}          \label{ss.e.reg}
The asymptotic upper bounds (\ref{e.e.low.conf.bd/cone})
and~(\ref{e.e.2sconf.bd/lin}) for the confidence probabilities derived in
Theorem~\ref{t:e:conf.bd} seem to involve only~$P$. The model~${\cal P}$
and its tangent set~${\cal G}$ at~$P$, however, enter through the
regularity condition. As indicated above, the bounds are not meaningful
without such regularity conditions.
\subsubsection{Modified Regularity, Asymptotic Linearity and Normality}
\paragraph{Asymptotic Median Bias}
In Theorem~\ref{t:e:conf.bd}, the regularity conditions
  (\ref{e.e.k.oomedS.ge0/cone}), (\ref{e.e.k.oomedS.ge0/lin})
and~(\ref{e.e.k.oomedS.le0/lin}), respectively, are certainly fulfilled
if asymptotic median nonnegativity, respectively nonpositivity, holds
for every fixed tangent in the respective tangent set~${\cal G}$,
in the sense that
\begingroup \mathsurround0em\arraycolsep0em \begin{eqnarray}
\label{e.e.oomedS.ge0.g/all}
   \Nliminf_{n} P_{n,t_n,g}^n \bigl\{
   S_n \ge T(P_{n,t_n,g})  \bigr\} &{}\ge{}& \Tfrac{1}{2} \\
\label{e.e.oomedS.le0.g/all}
   \Nliminf_{n} P_{n,t_n,g}^n \bigl\{
   S_n \le T(P_{n,t_n,g})  \bigr\} &{}\ge{}& \Tfrac{1}{2}
\end{eqnarray}\endgroup
respectively, for every~$g\in{\cal G}$ and every convergent sequence
  $ t_n\to t\in (0,\infty)$. The notion implicitly depends on~$P$,
the model~${\cal P}$, and its tangent set~${\cal G}$ at~$P$.
\par     Asymptotic median unbiasedness, that is,
         (\ref{e.e.oomedS.ge0.g/all}) and~(\ref{e.e.oomedS.le0.g/all}),
for every $g\in {\cal G}$, is the regularity assumption
of Pfanzagl and Wefelmeyer~(1982; Theorem~9.2.2).
\paragraph{Asymptotic Linear Estimators}
An estimator sequence~$(S_n)$ is asymptotically linear at~$P$ if there
exists some function $\eta \in L_2(P)\cap \{\mbox{const}\}^{\perp} $,
the (unique) influence curve of~$(S_n)$ at~$P$, such that
\begin{equation} \label{e.e.asLest.IC}
    \sqrt{n}\,\bigl(S_n-T(P)\bigr) =
    \rave_{i=1}^n \eta(x_i)
    + \Lo_{P^n}(n^0)
\end{equation}
For example, the estimator sequences $(\bar{S}_n)$ and~$(\tilde{S}_n)$,
in view of~(\ref{e.e.barS}) and~(\ref{e.e.tildeS}), are asymptotically
linear at~$P$ with influence curves~$\bar{\kappa}$ and~$\tilde{\kappa}$,
respectively.
\par
The construction of such estimators, given a family of prescribed
influence curves~$\eta_P$, 
one for each (unknown) element~$P$ of model~${\cal P}$, under very
general conditions, 
is a topic in Bickel et al.~(1993; Appendix~A.10),
van der Vaart~(1998; Sections~25.8--10), and HR~(1994; Chapters~1 and~6).
\subparagraph{Asymptotic Normality}
The expansions~(\ref{e.e.asLest.IC}), (\ref{e.e.sntkg}),
and~(\ref{E:e:llh.tn/+/-}),
of the estimator, the functional, and loglikelihoods, respectively, imply
the following asymptotic normality extending~(\ref{e.e.asyNhatS}),
\begin{equation} \label{e.e.asyNhatSgt}
   \Bigl \lgroup \sqrt{n}\,\bigl(S_n-T(P_{n,t_n,g})\bigr)
   \Bigr \rgroup   \bigl( P_{n,t_n,g}^n\bigr) \gwto {\cal N}\bigl(
   \hspace{3\ei} t \hspace{6\ei}\langle \eta-\kappa|g\rangle,
   \Vert \eta\Vert^2 \hspace{6\ei}\bigr)
\end{equation}
for all convergent $ t_n\to t\in (0,\infty)$, every $g\in {\cal G}$,
and so, for each $c\in [\hskip6\ei0,\infty)$,
\begingroup \mathsurround0em\arraycolsep0em \begin{eqnarray}
\hspace{-2.5em}  \label{e.e.Shat.Pgt<0}
  \Nlim_{n} P_{n,t_n,g}^n \bigl\{
  \sqrt{n}\,\bigl(S_n- T(P_{n,t_n,g})\bigr) < c \,\bigr\} &{}={}&
  \Phi \Bigl(\,\frac{c-t\langle \eta-\kappa|g\rangle}{\Vert
      \eta\Vert }\,\Bigr) \\
\hspace{-2.5em} \label{e.e.Shat.Pgt>0}
  \Nlim_{n} P_{n,t_n,g}^n \bigl\{
  \sqrt{n}\,\bigl(S_n- T(P_{n,t_n,g})\bigr) > -c \,\bigr\} &{}={}&
  \Phi \Bigl(\, \frac{c+t\langle \eta-\kappa|g\rangle}{\Vert
      \eta\Vert }\,\Bigr)
\end{eqnarray}\endgroup
where of course ${}< c$ may also be replaced by~${}\le c$.
These convergences in particular apply to $(\bar{S}_n)$ and~$(\tilde{S}_n)$,
with $\eta=\bar{\kappa}$, respectively~$\eta= \tilde{\kappa}$.
\subparagraph{Asymptotic Confidence Probabilities Based \boldmath
              on~$(\breve{S}_n)$:}
Besides~$(\tilde{S}_n)$, we consider any optimal estimator
sequences~$(\breve{S}_n)$ as described by~(\ref{e.e.u.Shat.cone}).
Then, by the asymptotic normality~(\ref{e.e.asyNhatSgt}) of~$(\tilde{S}_n)$,
and contiguity, we conclude that, for all convergent $t_n\to t$
in~$(0,\infty)$, every tangent~$g\in\tilde{{\cal G}}$, and each $c\in\R$,
\begingroup \mathsurround0em\arraycolsep0em \begin{eqnarray}
\lefteqn{\hskip-2.25em      P_{n,t_n,g}^n \bigl\{
     \sqrt{n}\,\bigl(\breve{S}_n- T(P_{n,t_n,g})\bigr) < c \,\bigr\} }
\label{e.cr.wcvSbrvStl0}\\      &{}={}&
     P_{n,t_n,g}^n \bigl\{ \sqrt{n}\,\bigl(\breve{S}_n- T(P)\bigr)
     <  c + t \hspace{6\ei} \langle \kappa|g \rangle + \Lo(n^0) \bigr\}
\hspace{-4.5em}\nonumber \\      &{}={}&
     P_{n,t_n,g}^n \bigl\{  \sqrt{n}\,\bigl(\breve{S}_n- T(P)\bigr)_+
     <  c + t \hspace{6\ei} \langle \kappa|g \rangle + \Lo(n^0) \bigr\}
\hspace{-4.5em}\nonumber \\      &{}={}&
     P_{n,t_n,g}^n \bigl\{  \sqrt{n}\,\bigl(\tilde{S}_n- T(P)\bigr)_+ <
     c + t \hspace{6\ei} \langle \kappa|g \rangle + \Lo_{P^n}(n^0) \bigr\}
\hspace{-4.5em}\nonumber \\      &{}={}&
     P_{n,t_n,g}^n \bigl\{  \sqrt{n}\,\bigl(\tilde{S}_n- T(P)\bigr)_+
     <  c + t \hspace{6\ei} \langle \kappa|g \rangle \bigr\} + \Lo(n^0)
\hspace{-4.5em}\nonumber \\       &{}={}&  
     P_{n,t_n,g}^n \bigl\{  \sqrt{n}\,\bigl(\tilde{S}_n- T(P)\bigr)
     <  c + t \hspace{6\ei} \langle \kappa|g \rangle \bigr\} + \Lo(n^0)
\hspace{-4.5em}\nonumber \\       &{}={}&  
     P_{n,t_n,g}^n \bigl\{
     \sqrt{n}\,\bigl(\tilde{S}_n- T(P_{n,t_n,g})\bigr) < c \,\bigr\}
     + \Lo(n^0) \label{e.cr.wcvSbrvStl01}
\end{eqnarray}\endgroup
provided that
\begin{equation} \label{e.cr.wcvSbrvStl02}
   t \hspace{6\ei} \langle \kappa|g \rangle > -c \hspace{-.5em}
\end{equation}
In~(\ref{e.cr.wcvSbrvStl0})--(\ref{e.cr.wcvSbrvStl01}),
we may replace $\ran{}<{}\ran$ 
            by~$\ran{}\le {}\ran$, 
hence by any inequality sign.
\par
Thus (\ref{e.e.Shat.Pgt<0}) and~(\ref{e.e.Shat.Pgt>0}),
using~$-c$ instead of~$c$, extend from~$(\tilde{S}_n)$ to~$(\breve{S}_n)$.
\subsubsection{One-Sided Regularity \boldmath
               of~$(\tilde{S}_n)$, $(\breve{S}_n)$}
\label{sss.cr.reg}
Let $c=0$, as in the regularity assumptions of Theorem~\ref{t:e:conf.bd}.
\paragraph{Regularity \boldmath
           of~$(\bar{S}_n)$:}
For $\eta= \bar{\kappa}$, since
    $ \langle \bar{\kappa}-\kappa|g\rangle=0
      \enskip \forall g\in \bar{{\cal G}} $,
the two limits in (\ref{e.e.Shat.Pgt<0}) and~(\ref{e.e.Shat.Pgt>0})
are always~$1/2$. Hence $(\bar{S}_n)$ is asymptotically median unbiased,
for each $g\in \bar{{\cal G}}$. In particular, $(\bar{S}_n)$ satisfies
conditions~(\ref{e.e.k.oomedS.ge0/lin}) and~(\ref{e.e.k.oomedS.le0/lin}).
\paragraph{Regularity \boldmath
           of~$(\tilde{S}_n)$, $(\breve{S}_n)$:}
For $\eta=\tilde{\kappa}$, since $t>0$ and
    $ \langle \tilde{\kappa}-\kappa|g\rangle \ge0 \enskip
      \forall g\in \tilde{{\cal G}} $,
the limit in~(\ref{e.e.Shat.Pgt<0}) is always~$\le 1/2$
(it is $=1/2$, e.g.\ for\footnote{In Subsections~\ref{ss.e.reg}
  and~\ref{ss.e.lbh}, the choice $g=\tilde{\kappa}$ stands under
  the provision that $\tilde{\kappa}\in \tilde{{\cal G}}$.}
  $g=0,\tilde{\kappa}$).
\par     Thus, $(\tilde{S}_n)$ satisfies the asymptotic median nonnegativity
condition~(\ref{e.e.oomedS.ge0.g/all}), for every~$g\in \tilde{{\cal G}}$,
hence, in particular,
         $(\tilde{S}_n)$ fulfills condition~(\ref{e.e.k.oomedS.ge0/cone}).
\par     If~$(\breve{S}_n)$ satisfying~(\ref{e.e.u.Shat.cone}) is another
optimal estimator sequence,
  (\ref{e.cr.wcvSbrvStl0})--(\ref{e.cr.wcvSbrvStl02}) apply with $c=0$,
hence asymptotic median nonnegativity~(\ref{e.e.oomedS.ge0.g/all})
of~$(\tilde{S}_n)$ is inherited to~$(\breve{S}_n)$,
for every~$g\in \nolinebreak\tilde{{\cal G}}$ such that
  $ \langle \kappa|g\rangle>0$. This suffices to fulfill
  condition~(\ref{e.e.k.oomedS.ge0/cone}), because
  $ \langle \kappa|\tilde{\kappa}\rangle = \Vert \tilde{\kappa}\Vert^2>0 $,
  and so eventually $ \langle \kappa|g_m\rangle >0 $ for any tangents
  $ g_m\in \tilde{{\cal G}}$ approaching~$\tilde{\kappa}$.
\subsubsection{Positive Median Bias \boldmath
               of~$(\tilde{S}_n)$, $(\breve{S}_n)$}
\label{sss.cr.posmedbias}
For $c=0$ and $\eta=\tilde{\kappa}$, the limit in~(\ref{e.e.Shat.Pgt>0})
($=1/2$ for $g=0,\tilde{\kappa}$) in general falls below~$1/2$. We shall
prove this for any estimator sequence~$(\breve{S}_n)$ which is optimal
in the sense of Theorem~\ref{t:e:conf.bd}(\ref{i.e.cone}).
\par
Consequently, all these estimators violate asymptotic median
nonpositivity~(\ref{e.e.oomedS.le0.g/all}). The result corresponds to
the level breakdown we have encountered in Subsection~\ref{ss.t.Lbruch}.
\begin{Prop}\sl \label{p.cr.posmed}
Let~$\tilde{{\cal G}}$ be a convex tangent cone such that
\begin{equation} \label{e.ccec.bknetlk}
    \bar{\kappa}\ne \tilde{\kappa}
\end{equation}
Then there is some tangent $g_0\in \tilde{{\cal G}}$ such that
  $ \langle \kappa|g_0\rangle > 0 $ and
\begin{equation}\label{e.e.asybias=1/cone}
  \Ninf_{t>0} \Nlim_{n} P_{n,t,g_0}^n \bigl\{
  \breve{S}_n \le T(P_{n,t,g_0})   \bigr\} = 0
\end{equation}
for all estimator sequences~$(\breve{S}_n)$
of kind\/~{\rm (\ref{e.e.u.Shat.cone}).}
\end{Prop}
\begin{Bew}\enskip
If $\bar{\kappa}\ne \tilde{\kappa}$ there is some tangent
   $g_0\in \tilde{{\cal G}}$ such that
\begin{equation} \label{e.e.kgtlkg}
  \langle \kappa|g_0 \rangle < \langle \tilde{\kappa}|g_0 \rangle
\end{equation}
Then, for $c=0$ and $\eta=\tilde{\kappa}$, (\ref{e.e.Shat.Pgt>0}) implies
\begin{equation}\label{e.e.asybias=1/coneStl}
  \Ninf_{t>0} \Nlim_{n} P_{n,t,g_0}^n \bigl\{
  \tilde{S}_n \le T(P_{n,t,g_0})   \bigr\}
  = \inf_{t>0}\:
    \Phi \Bigl(-\frac{t\langle \tilde{\kappa}-\kappa|g\rangle}{\Vert
    \tilde{\kappa}\Vert }\,\Bigr) = 0
\end{equation}
But~$g_0$ may always be chosen such that, in addition to~(\ref{e.e.kgtlkg}),
\begin{equation}\label{e.e.0<kgtlkg}
   0< \langle \kappa|g_0 \rangle < \langle \tilde{\kappa}|g_0 \rangle
\end{equation}
If necessary, pass to a suitable convex combination~$g_{02}$ of~$g_0$
satisfying~(\ref{e.e.kgtlkg}) and~$\tilde{\kappa}$, in order to
achieve~(\ref{e.e.0<kgtlkg}).
\par     Then the arguments
         (\ref{e.cr.wcvSbrvStl0})--(\ref{e.cr.wcvSbrvStl02})
go through, with $c=0$, and with~${}\le{}$ in the place of~${}<{}$.
Thus the positive asymptotic median bias~(\ref{e.e.asybias=1/cone})
carries over from~$(\tilde{S}_n)$ to all estimator
sequences~$(\breve{S}_n)$ satisfying~(\ref{e.e.u.Shat.cone}).\qed
\end{Bew}
\begin{Rem}\rm  \label{r.cr.pfanz}
The result implies that bound~(\ref{e.e.low.conf.bd/cone}) cannot
possibly be achieved if, in addition to~(\ref{e.e.k.oomedS.ge0/cone}),
asymptotic median nonpositivity~(\ref{e.e.oomedS.le0.g/all}) is imposed
(for all~$g\in \tilde{{\cal G}}$, or only all $g\in \tilde{{\cal G}}$
such that $\langle \kappa|g \rangle>0$).
In particular, asymptotic median unbiasedness cannot be afforded if
bound~(\ref{e.e.low.conf.bd/cone}) is to be attained.
\par
As a consequence, Theorem~9.2.2 of Pfanzagl and Wefelmeyer~(1982)
for (closed) convex tangent cones~$\tilde{{\cal G}}$ is ailing
in two respects:
\par
First, since $-g\notin \tilde{{\cal G}}$ in general and
             $-\tilde{\kappa}\notin \tilde{{\cal G}}$ in particular,
their bound of form~(\ref{e.e.2sconf.bd/lin}) for two-sided confidence
limits, with~$\tilde{\kappa}$ in the place of~$\bar{\kappa}$, is not
available over cones, but only bound~(\ref{e.e.low.conf.bd/cone}) for
lower confidence limits.
\par
Second, their regularity condition is too strict:
Contrary to what they believe (Section~9.1, p~154), even the one-sided
bound~(\ref{e.e.low.conf.bd/cone}), let alone the asserted two-sided
extension~(\ref{e.e.2sconf.bd/lin}), cannot possibly be achieved by any
estimator sequence that is asymptotically median unbiased.\qed \end{Rem}
\subsection{Comparison of Cones and Spaces}               \label{ss.e.vs}
\paragraph{Variance and Sample Size}
Recall the setup of Subsection~\ref{ss.t.pow.cone/lin}:
  $P\in \tilde{{\cal P}}\subset \bar{{\cal P}}$, with tangent
set a convex cone~$\tilde{{\cal G}}$, respectively the linear span
(\ref{e.t.Gbar=spanGtil}):~$ \bar{{\cal G}}=\lin \tilde{{\cal G}}$.
\par
Then, in view of the asymptotic normality~(\ref{e.e.asyNhatS}),
the previous comparison of 
  $\Vert \tilde{\kappa}\Vert^2$ and~$\Vert \tilde{\kappa}\Vert^2$
now concerns the variances~$\Vert \tilde{\kappa}\Vert^2\!\big/n$
and~$\Vert \bar{\kappa}\Vert^2\!\big/n$ of the approximate normal
distributions of~$\tilde{S}_n-T(P)$ and~$\bar{S}_n-T(P)$, respectively.
\par
Thus, the value~$T(P)$, in terms of variance or width of confidence
intervals, can be estimated under~$P^n$ more accurately in model~${\cal P}$
with tangent set~$\tilde{{\cal G}}$ than it is possible in the larger
model~$\bar{{\cal P}}$ with tangent set~$\bar{{\cal G}}$.
Observations at the higher rate
  $  \bar{n}/\tilde{n} \to  \Vert \bar{\kappa}\Vert^2\!
                            \big/\Vert \tilde{\kappa}\Vert^2  $
are needed under~$P$ to estimate~$T(P)$ with the same asymptotic accuracy
    by~$\bar{S}_{\bar{n}}$ as by~$\tilde{S}_{\tilde{n}}$.
Again Example~\ref{ex1.t.kbkt/cone} applies.
\paragraph{Lower Confidence Limits for Spaces}
The preceding comparison does not explicitly take the different
sets of regularity assumptions into account:
In the case of~$(\tilde{S}_n)$,
it is condition~(\ref{e.e.k.oomedS.ge0/cone}), and conditions
  (\ref{e.e.k.oomedS.ge0/lin}),~(\ref{e.e.k.oomedS.le0/lin})
   in the case of~$(\bar{S}_n)$.
\par
However, in the case of a linear tangent space~$\bar{{\cal G}}$,
suppose we dispense of condition~(\ref{e.e.k.oomedS.le0/lin}) and,
keeping~(\ref{e.e.k.oomedS.ge0/lin}), wish to maximize the asymptotic
confidence probability merely of the sequence of lower confidence limits
$S_n-c/\!\sqrt{n}\,$ of~$T(P)$, under~$(P^n)$.
In particular, the statistical task seems to be made easier.
\par    Nevertheless, the previous upper bound
        $\Phi \bigl(c/\Vert \bar{\kappa}\Vert \hspace{6\ei}\bigr)$
established under Theorem~\ref{t:e:conf.bd}(\ref{i.e.lin}),
with \mbox{$t_n'=t'=\infty$,} does not increase,
and~$(\bar{S}_n)$ remains an optimal estimator sequence.
This is true, simply because~$\bar{{\cal G}}$ also is a convex tangent
cone, so Theorem~\ref{t:e:conf.bd}(\ref{i.e.cone}) especially
holds for~$\bar{{\cal G}}$.
\par
Thus, under 
      condition~(\ref{e.e.k.oomedS.ge0/lin}), asymptotic median
nonpositivity~(\ref{e.e.oomedS.le0.g/all}) for~$\bar{{\cal G}}$,
as well as the maximization, subject to~(\ref{e.e.k.oomedS.le0/lin}),
of the asymptotic confidence probability under~$(P^n)$ of the sequence
of upper confidence limits $S_n+c/\!\sqrt{n}\,$ for~$T(P)$, come free
with~$(\bar{S}_n)$, which achieves the corresponding upper bound, which
is~$\Phi \bigl(c/\Vert \bar{\kappa}\Vert \hspace{6\ei}\bigr)$ again.
\paragraph{Two-Sided Confidence Limits for Cones}
In the case of a convex tangent cone~$\tilde{{\cal G}}$, suppose we want
to maximize the asymptotic confidence probability under~$(P^n)$ of the
sequence of lower confidence limits $S_n-c/\!\sqrt{n}\,$ of~$T(P)$,
as in Theorem~\ref{t:e:conf.bd}(\ref{i.e.cone}), but insist on asymptotic
median unbiasedness, that is, (\ref{e.e.oomedS.ge0.g/all})
and~(\ref{e.e.oomedS.le0.g/all}) for every $g\in\tilde{{\cal G}}$.
As (\ref{e.e.oomedS.ge0.g/all}) and~(\ref{e.e.oomedS.le0.g/all})
imply~(\ref{e.e.k.oomedS.ge0/cone}), the statistical task is made
more difficult, and one expects the upper
bound~$\Phi \bigl(c/\Vert \tilde{\kappa}\Vert \,\bigr)$ to decrease.
According to Proposition~\ref{p.cr.posmed}, it must strictly decrease
if $\bar{\kappa}\ne\tilde{\kappa}$.
\par
We clarify the amount of decrease, at least in the class of estimator
sequences~$(S_n)$ which are asymptotically linear at~$P$.
For such an estimator with influence curve~$\eta$ at~$P$,
the lower/upper confidence limits $ S_n \mp c/\!\sqrt{n}\, $ 
satisfy
\begin{equation} \label{e.e.lim.confP.c.2}
\mathsurround0em\arraycolsep0em \begin{array}{r}\Ds
     \Nlim_{n} P^n \bigl\{ \sqrt{n}\,\bigl( S_n-T(P)\bigr)
     < c \,\bigr\} = \Phi \Bigl(\frac{c}{\Vert \eta\Vert }\Bigr)  \\
\Ds \rule{0pt}{4ex} {} =
     \Nlim_{n} P^n \bigl\{ \sqrt{n}\,\bigl( S_n-T(P)\bigr)
     > -c \,\bigr\} \end{array}
\end{equation}
Under local alternatives, in~view of the limits (\ref{e.e.Shat.Pgt<0})
and~(\ref{e.e.Shat.Pgt>0}) for each $g\in \tilde{{\cal G}}$, $(S_n)$~is
asymptotically median unbiased iff
  $  \langle \eta|g\rangle = \langle \kappa|g\rangle \enskip
        \forall g\in \tilde{{\cal G}}  $, which holds if and only if
\begin{equation}\label{e.e.sprodk.1/clincone}
   \langle \eta|g\rangle = \langle \kappa|g\rangle
   \qquad  \forall g\in \clin\tilde{{\cal G}}  \hspace{-1.5em}
\end{equation}
Introducing the projections $\bar{\kappa}$ of~$\kappa$,
and~$\bar{\eta}$ of~$\eta$, on~$ \clin \tilde{{\cal G}}$,
$\bar{\eta}$~must equal~$\bar{\kappa}$.
But, subject to $\bar{\eta}=\bar{\kappa}$, the asymptotic confidence
     probability~$\Phi \bigl(c/\Vert \eta\Vert \,\bigr)$
is maximized iff~$\Vert \eta\Vert $ is minimized, which is the case
iff $\eta=\bar{\kappa}$.
\par     Therefore, in the class of estimator sequences which are
asymptotically linear at~$P$, the unique solution is the estimator
sequence~$(\bar{S}_n)$ with influence curve~$\bar{\kappa}$.
And the achievable upper bound decreases
from~$\Phi \bigl(c/\Vert \tilde{\kappa}\Vert \,\bigr)$
  to~$\Phi \bigl(c/\Vert \bar{\kappa}\Vert \,\bigr)$.
\par
So the answer to the corresponding (open) question raised for testing
in Remark~\ref{r.t.extH.test/cone} 
turns out negative in the estimation context. 
\par
In addition, in view of~(\ref{e.e.lim.confP.c.2}), the upper confidence
limits $ \bar{S}_n+c/\!\sqrt{n}\,$ of~$T(P)$ supplied by~$(\bar{S}_n)$
have the same asymptotic confidence
  probability~$\Phi \bigl(c/\Vert \bar{\kappa}\Vert \,\bigr)$
under~$P^n$ as the lower confidence limits $ \bar{S}_n-c/\!\sqrt{n}\,$.
And the two-sided bounds~$ \bar{S}_n\mp c/\!\sqrt{n}\,$,
in view of~(\ref{e.e.lim.confPtg.c.2/lin/Sbar}) below, maintain their
asymptotic confidence probability for~$T$ even under local
perturbations~$P_{n,t,g}$ of~$P$, $g\in \tilde{{\cal G}}$.
\subsection{Local Behaviour of Efficient Confidence Limits} \label{ss.e.lbh}
\subsubsection{Confidence Probabilities Under Perturbations}
Given $c\in (0,\infty)$, we study the two sequences of lower/upper limits
$\bar{S}_n\mp c/\!\sqrt{n}\,$ and $\tilde{S}_n\mp c/\!\sqrt{n}\,$
under local perturbations~$P_{n,t,g}$ of~$P$.
\paragraph{Stability of Confidence Limits Based \boldmath
           on~$(\bar{S}_n)$:}
For $\eta= \bar{\kappa}$,
the two limits in (\ref{e.e.Shat.Pgt<0}) and~(\ref{e.e.Shat.Pgt>0}),
since $ \langle \bar{\kappa}-\kappa|g\rangle=0
        \enskip \forall g\in \bar{{\cal G}} $, are always the same,
\begin{equation} \label{e.e.lim.confPtg.c.2/lin/Sbar}
\mathsurround0em\arraycolsep0em \begin{array}{r}
\Ds   \Nlim_{n} P_{n,t_n,g}^n \bigl\{
      \sqrt{n}\,\bigl(\bar{S}_n-T(P_{n,t_n,g})\bigr) < c \,\bigr\} =
      \Phi \Bigl(\frac{c}{\Vert \bar{\kappa}\Vert }\Bigr) \\
\Ds   \rule{0pt}{4ex} {}=
      \Nlim_{n} P_{n,t_n,g}^n \bigl\{
      \sqrt{n}\,\bigl(\bar{S}_n-T(P_{n,t_n,g})\bigr) > -c \,\bigr\}
\end{array}\end{equation}
for every convergent sequence $ t_n\to t$ in~$(0,\infty)$, every
$g\in \bar{{\cal G}}$, which reveals some stability of the lower/upper
limits based on~$(\bar{S}_n)$.
\paragraph{Instability of Confidence Limits Based \boldmath
           on~$(\tilde{S}_n)$, $(\breve{S}_n)$:}
For $\eta=\tilde{\kappa}$, the limits in~(\ref{e.e.Shat.Pgt<0})
and~(\ref{e.e.Shat.Pgt>0}) are, respectively,
\begingroup \mathsurround0em\arraycolsep0em \begin{eqnarray}
\hspace{-2.5em}  \label{e.e.Stild.Pgt<0}
   \Nlim_{n} P_{n,t_n,g}^n \bigl\{
   \sqrt{n}\,\bigl(\tilde{S}_n-T(P_{n,t_n,g})\bigr) < c \,\bigr\} &{}={}&
   \Phi \Bigl(\, \frac{c-t\langle \tilde{\kappa}-\kappa|g\rangle}{\Vert
          \tilde{\kappa}\Vert }\,\Bigr)  \\
\hspace{-2.5em} \label{e.e.Stild.Pgt>0}
  \Nlim_{n} P_{n,t_n,g}^n \bigl\{
  \sqrt{n}\,\bigl(\tilde{S}_n- T(P_{n,t_n,g})\bigr) > -c \,\bigr\} &{}={}&
  \Phi \Bigl(\, \frac{c+t\langle \tilde{\kappa}-\kappa|g\rangle}{\Vert
          \tilde{\kappa}\Vert }\,\Bigr)
\end{eqnarray}\endgroup
\subparagraph{Under-Coverage by Lower Confidence Limits}
The limit~(\ref{e.e.Stild.Pgt<0}) is always
    $ {}\le \Phi \bigl(c/\Vert \tilde{\kappa}\Vert \hspace{6\ei}\bigr) $,
since $t>0$ and $ \langle \tilde{\kappa}-\kappa|g\rangle \ge0 \enskip
                  \forall  g\in \tilde{{\cal G}} $;
the upper bound is achieved, e.g.~for $g=0,\tilde{\kappa}$.
In~general, e.g.\ for~$g_0$ taken from~(\ref{e.e.0<kgtlkg}), the limit
in~(\ref{e.e.Stild.Pgt<0}), with~${}\le c$ in the place of~${}<c$,
may become arbitrarily close to~$0$ as
\begin{equation} \label{e.e.asybias.Stild=-1.c/cone}
  \Ninf_{t>0} \Nlim_{n} P_{n,t,g_0}^n \bigl\{
  \sqrt{n}\,\bigl(\tilde{S}_n-T(P_{n,t,g_0})\bigr) \le c \,\bigr\} = 0
\end{equation}
In view of (\ref{e.cr.wcvSbrvStl0})--(\ref{e.cr.wcvSbrvStl02}),
the limit statement~(\ref{e.e.Stild.Pgt<0}) for
  $ t\hspace{4\ei}\langle \kappa|g \rangle >-c$ extends to~$(\breve{S}_n)$,
hence also~(\ref{e.e.asybias.Stild=-1.c/cone}) extends to~$(\breve{S}_n)$.
Obviously, (\ref{e.e.asybias.Stild=-1.c/cone})
           generalizes~(\ref{e.e.asybias=1/cone}).
\subparagraph{Over-Coverage by Upper Confidence Limits}
The limit in~(\ref{e.e.Stild.Pgt>0}) is
always~$\ge \Phi \bigl(c/\Vert \tilde{\kappa}\Vert \hspace{6\ei}\bigr)$;
and~${}=\Phi \bigl(c/\Vert \tilde{\kappa}\Vert \hspace{6\ei}\bigr)$
e.g.\ for~$g=0,\tilde{\kappa}$. In~general, e.g.\ for~$g_0$ taken
from~(\ref{e.e.0<kgtlkg}), the limit in~(\ref{e.e.Stild.Pgt>0})
may become arbitrarily close to~$1$,
\begin{equation} \label{e.e.asybias.Stild=1.-c/cone}
  \Nsup_{t>0}   \Nlim_{n} P_{n,t,g_0}^n \bigl\{
  \sqrt{n}\,\bigl(\tilde{S}_n- T(P_{n,t,g_0})\bigr) > -c \,\bigr\} = 1
\end{equation}
In view of (\ref{e.cr.wcvSbrvStl0})--(\ref{e.cr.wcvSbrvStl02}),
with $-c$ in the place of~$c$, the limit statement~(\ref{e.e.Stild.Pgt>0})
extends from~$(\tilde{S}_n)$ to~$(\breve{S}_n)$ of
form~(\ref{e.e.u.Shat.cone}), provided that
   $ t \hspace{6\ei}\langle \kappa|g_0\rangle>c $, and hence
also~(\ref{e.e.asybias.Stild=1.-c/cone}) extends to~$(\breve{S}_n)$.
\par
The degenerate limits (\ref{e.e.asybias.Stild=-1.c/cone})
and~(\ref{e.e.asybias.Stild=1.-c/cone}) indicate an instability
of the lower and upper confidence limits based on~$(\tilde{S}_n)$,
which is not accounted for by the criterion maximized in
Theorem~\ref{t:e:conf.bd}(\ref{i.e.cone}) merely under~$(P^n)$,
nor by the (only one-sided) asymptotic median nonnegativity conditions
(\ref{e.e.k.oomedS.ge0/cone}) or~(\ref{e.e.oomedS.ge0.g/all}).
\subsubsection{\ran \boldmath
                $(\bar{S}_n)$, $(\tilde{S}_n)$, $(\breve{S}_n)$
                in the Light of the Convolution Theorem}
\paragraph{Superefficiency}
The convolution theorem by van der Vaart~(1998; Theorem~25.20) states the
lower bound~$ \Vert \bar{\kappa}\Vert^2 $ for the asymptotic variance, which
is attained by~$(\bar{S}_n)$ in~(\ref{e.e.asyNhatS}), but which seems to
contradict the smaller asymptotic variance~$ \Vert \tilde{\kappa}\Vert^2 $
of~$(\tilde{S}_n)$ in~(\ref{e.e.asyNhatS}),
in case~(\ref{e.ccec.bknetlk}): $\bar{\kappa}\ne \tilde{\kappa}$.
\paragraph{H\'ajek--Regularity}
This convolution result concerns the asymptotic variance of estimator
sequences~$(S_n)$ which are H\'ajek--regular.
$(S_n)$~is called H\'ajek--regular at~$P$, for the functional~$T$, along the
tangent set~${\cal G}$, if there is some (limit) distribution~$M$ such that,
for every $g\in \tilde{{\cal G}}$ and $ t_n\to t$ in~$(0,\infty)$,
\begin{equation} \label{e.e.asyMSgt}
   \Bigl \lgroup \sqrt{n}\,\bigl(S_n-T(P_{n,t_n,g})\bigr) \Bigr \rgroup
   \bigl( P_{n,t_n,g}^n\bigr) \gwto M
\end{equation}
If~$(S_n)$ is H\'ajek--regular with limit~$M$, then $M(0,\infty)\ge 1/2$
implies asymptotic median nonnegativity~(\ref{e.e.oomedS.ge0.g/all}),
   $M(-\infty,0)\ge 1/2$ implies asymptotic median
   nonpositivity~(\ref{e.e.oomedS.le0.g/all}),
and $M(0,\infty)= 1/2$, $M(\{0\})=0$, implies
that $(S_n)$ is asymptotically median unbiased.
\paragraph{H\'ajek--Nonregularity}
Contrary to~$(\bar{S}_n)$, whose limit distribution in~(\ref{e.e.asyNhatSgt})
is always~${\cal N}\bigl(0,\Vert \bar{\kappa}\Vert^2 \hspace{6\ei}\bigr)$,
hence is H\'ajek--regular, the limit distribution of~$(\tilde{S}_n)$
in~(\ref{e.e.asyNhatSgt}) clearly does depend on the
        particular~$(t,g)\in (0,\infty)\times\tilde{{\cal G}}$.
Therefore, the estimator sequence~$(\tilde{S}_n)$ is not H\'ajek--regular.
As (\ref{e.e.Stild.Pgt<0}),~(\ref{e.e.asybias.Stild=-1.c/cone}) with
  $c\in (0,\infty)$ also hold for~$(\breve{S}_n)$ and $g=0$, respectively
for the tangent~$g_0$ taken from~(\ref{e.e.0<kgtlkg}), neither estimator
sequence~$(\breve{S}_n)$ which is optimal in the sense of
  Theorem~\ref{t:e:conf.bd}(\ref{i.e.cone}) can be H\'ajek--regular,
unless $ \tilde{\kappa}= \bar{\kappa}$.
\section{Appendix}                    \setcounter{equation}{0}\label{s.X}
\subsection{Projection---Generalities}
Let ${\cal H}$ be a Hilbert space---for example, ${\cal H}=L_2(P)$---and
fix some $\kappa\in {\cal H}$.
\par
If $\bar{{\cal G}}$ is a closed linear subspace of~${\cal H}$,
the orthogonal projection~$\bar{\kappa}\in \bar{{\cal G}}$
of~$\kappa$ on~$\bar{{\cal G}}$, and unique element of~$\bar{{\cal G}}$
closest to~$\kappa$ in norm~$\Vert \makebox[.5em][c]{.}\Vert $,
is characterized by
\begin{equation} \label{e:X:lin}
   \langle \kappa - \bar{\kappa}|g \rangle = 0
   \quad \forall g\in \bar{{\cal G}}
\end{equation}
If~$\tilde{{\cal G}}$ is a closed convex cone in~${\cal H}$,
the projection~$\tilde{\kappa}\in \tilde{{\cal G}}$ of~$\kappa$
on~$\tilde{{\cal G}}$, that is, the unique element of~$\tilde{{\cal G}}$
closest to~$\kappa$ in norm~$\Vert \makebox[.5em][c]{.}\Vert $, is
characterized by
\begin{equation} \label{e:X:cone}
   \langle \kappa |\tilde{\kappa} \rangle =
   {\Vert \tilde{\kappa}\Vert}^2\weg,    \qquad
   \langle \kappa |g \rangle \le \langle \tilde{\kappa} |g \rangle
   \quad \forall g\in \tilde{{\cal G}}
\end{equation}
If~$\hat{\cal G} $ is an arbitrary nonempty closed convex subset
of~${\cal H}$, the 
unique minimum norm element~$\hat{g}$ of~$\hat{{\cal G}}$
is characterized by
\begin{equation} \label{e:X:cvx}
     {\Vert \hat{g}\Vert}^2 \le  \langle g|\hat{g} \rangle
    \quad \forall g\in \hat{\cal G}
\end{equation}
These facts are well-known; see, for example, Proposition~4.2.1
in Pfanzagl and Wefelmeyer (1982).
(\ref{e:X:cvx})~may be proved by differentiation at $s=0$ of the function
  $ {\Vert (1-s)\hat{g}+sg \Vert}^2 $, which is convex in $ 0\le s\le1 $,
for any $g\in \nolinebreak\hat{\cal G}$.
Passing to~$\kappa- \tilde{{\cal G}}$ and using the structure of cones,
(\ref{e:X:cone})~may be derived from~(\ref{e:X:cvx}).
Using~$ - \bar{{\cal G}}= \bar{{\cal G}}$ for the linear
space~$\bar{{\cal G}}$, (\ref{e:X:lin}) follows from~(\ref{e:X:cone}).
\subsection{Projection---Examples} 
\paragraph{{\rm ad} Example~\ref{ex1.t.kbkt/cone}:}
Recall (\ref{e.t.ex1.coeff}) and~(\ref{e.t.ex1.kb}). Then
\begin{equation}
  \bar{\gamma}_1>0 \iff b_1>b_2 \hspace{6\ei}c
  \iff \varphi(0)-\varphi(a)<\varphi(0)
\end{equation}
Introduce the function
  $ r(a)= \bigl[\varphi(0)-\varphi(a)\bigr]\big/
          \bigl[2 \hspace{6\ei}\Phi(a)-1\bigr]  $. Then
\begin{equation}
  \bar{\gamma}_2<0 \iff b_2 < b_1 c \iff r(a)< \varphi(0)
\end{equation}
However, $\varphi(0)=\lim_{a\to \infty}r(a)$ and
$\lim_{a \downarrow 0}r(a)=0$ (de l'Hospital). Moreover,
\begingroup \mathsurround0em\arraycolsep0em
\begin{eqnarray}
  & \Ds   \dot r(a)>0 \iff \varphi(0)-\varphi(a)<
  a \bigl[ \Phi(a)-\Tfrac{1}{2}\hspace{3\ei}\bigr] & \\
\noalign{\noindent But\nopagebreak} \label{e.x.dfi<adgfi}
  & \Ds   \varphi(0)-\varphi(a)=\itg_0^a x \hspace{6\ei}\varphi(x)\,dx
  < a \itg_0^a \varphi(x)\,dx
  = a \bigl[ \Phi(a)-\Tfrac{1}{2}\hspace{3\ei}\bigr] &
\end{eqnarray}\endgroup
Also $b_2<b_1$, since $b_2<b_1 c$ and $c<1$ (Cauchy--Schwarz).
\paragraph{{\rm ad} Example~\ref{ex2.t.JnotH/cone}:}
Recall $ b_1= \langle \kappa|g_1 \rangle = 2 \hspace{6\ei}\varphi(0) $
from~(\ref{e.t.ex1.coeff}), $g_3$ from~(\ref{e.t.ex2.g3}), and put
  $ b_3 = \langle \kappa|g_3 \rangle $, $ c= \langle g_1|g_3\rangle $.
Set $\sigma=\delta/\eta$. 
Then
\begin{equation} \label{e.x.g3N1}
   1 = {\Vert g_3\Vert}^2 = 2 \hspace{6\ei} \eta^2 \bigl(
   \sigma^2 \bigl[ \Phi(a)-\Tfrac{1}{2}\hspace{3\ei}\bigr]
   + \bigl[ 1-\Phi(a)\bigr]\bigr)
\end{equation}
As $ \frac{1}{2}b_3 = \delta \bigl[\varphi(0)-\varphi(a)\bigr]
                      - \eta \hspace{6\ei}\varphi(a) $, we have
\begingroup \mathsurround0em\arraycolsep0em \begin{eqnarray}
\label{e.x.sklein}  b_3<0 & {}\iff {}& \sigma \hspace{6\ei}
              \bigl[\varphi(0)-\varphi(a)\bigr] < \varphi(a) \\
\noalign{\vspace{\belowdisplayshortskip}\pagebreak[0]
\mathsurround=\msu\noindent And as $ \frac{1}{2}c =
          \delta \bigl[ \Phi(a)-\Tfrac{1}{2}\hspace{3\ei}\bigr]
           - \eta \bigl[ 1-\Phi(a)\bigr] $,
          we have\nopagebreak\vspace{\abovedisplayskip}}
\label{e.x.sgross}     c>0 & {}\iff{} &  \sigma \hspace{6\ei} \bigl[
    \Phi(a)-\Tfrac{1}{2}\hspace{3\ei}\bigr] > \bigl[ 1-\Phi(a)\bigr]
\end{eqnarray}\endgroup
But     
\begingroup \mathsurround0em\arraycolsep0em \begin{eqnarray}
\noalign{\vspace{-\abovedisplayskip}\vspace{\abovedisplayshortskip}}
\label{e.x.dfi>adgfi}  a \bigl[ 1-\Phi(a)\bigr]   &{}<{}&
    \itg_a^{\infty} x \hspace{6\ei}\varphi(x)\,dx =\varphi(a) \\
\noalign{\noindent and~(\ref{e.x.dfi<adgfi}),\nopagebreak}
  \frac{\varphi(0)-\varphi(a)}{\Phi(a)-\Phi(0)} &{}<{}& a <
  \frac{\varphi(a)-\varphi(\infty)}{\Phi(\infty)- \Phi(a)} \\
\noalign{\noindent imply\nopagebreak}
  \frac{\varphi(a)}{\varphi(0)-\varphi(a)} &{}>{}& \sigma  >
  \frac{1- \Phi(a)}{\Phi(a)- \Phi(0)}
\end{eqnarray}\endgroup 
for $ \sigma=\sigma_a=a \bigl[ 1-\Phi(a)\bigr]\big/
             \bigl[\varphi(0)-\varphi(a)\bigr] $.
Then~(\ref{e.x.g3N1}) defines us~$\eta=\eta_a$. 
\par     As $b_3<0<c,b_1$, the coefficients of the projection~$\bar{\kappa}$
on~$\bar{{\cal G}}=(\hspace{-15\ei}\cl)\lin \{g_1,g_2\}$ satisfy
   $ \bar{\gamma}_1>0>\bar{\gamma}_3$; confer~(\ref{e.t.ex1.kb}).
Therefore, $\bar{\kappa}\ne$~the projection~$\tilde{\kappa}$ on
the (closed) convex cone~$\tilde{{\cal G}}$ generated by~$g_1$ and~$g_3$,
and so $\Vert \tilde{\kappa}\Vert < \Vert \bar{\kappa}\Vert $.
\par     In minimizing the Lagrangian corresponding to~(\ref{e.t.ex1.Lagr}),
we can again rule out that both multipliers vanish. If $\beta_1>0$ then
  $\tilde{\gamma}_1=0$ and $\tilde{\gamma}_3=b_3+\beta_3\ge0$. As $b_3<0$,
necessarily $\beta_3>0$, hence $\tilde{\gamma}_3=0$ as
  $\beta_3\tilde{\gamma}_3=0$, which leads to an approximation error
  of~$\Vert \kappa-0\Vert^2= 1$.
This is worse than the error obtained under the assumption that $\beta_3>0$.
For in this case, $\tilde{\gamma}_3=0$ and $\tilde{\gamma}_1=b_1+\beta_1$
where $\beta_1=0$ due to $\beta_1 \tilde{\gamma}_1=0$ and $b_1>0$.
Hence $\tilde{\gamma}_1=b_1$, and the error amounts to
  $\Vert \kappa-b_1g_1\Vert^2=1-b_1^2 <1 $. Altogether,
this proves that $\tilde{\kappa}=b_1g_1$.
\subsection{Approximate Uniqueness} 
Given two probabilites $P$ and~$Q$ on some sample space,
let~$\tau^*$ be a Neyman--Pearson test for $P$ vs.~$Q$,
with critical value $c\in [\,0,\infty]$,
\begin{equation}
    \tau^* = \Jc (dQ>c \hspace{9\ei}dP \hspace{9\ei})
    \hspace{2em}\mbox{on}\hspace{.5em}
                   \{dQ\ne c \hspace{9\ei}dP \hspace{6\ei}\}
    \hspace{-2em}
\end{equation}
and possibly nonconstant randomization on
                  $\{dQ=c \hspace{9\ei}dP\hspace{6\ei}\}$.
By~$|\nu_c|$ we denote the total variation measure
of $d \nu_c=dQ-c \hspace{9\ei}dP$.
\begin{Lem}\sl \label{l.x.u}
Consider any test~$\tau$ for $P$ vs.~$Q$ such that,
for some {\rm $\delta \in (0,1)$,}
\begingroup \mathsurround0em\arraycolsep0em
\begin{eqnarray}  \label{e.x.lpt*t}
  & \Ds   \itg \tau\, dP\le \itg \tau^*\,dP+ \delta\weg, \qquad
          \itg \tau\, dQ\ge \itg \tau^*\,dQ - \delta    &  \\
\noalign{\noindent Then\nopagebreak}  \label{e.x.nct*t} & \Ds
  |\nu_c| \bigl\{ |\tau- \tau^*|>\varepsilon\bigr\}\le
  (1+c)\hspace{1\ei} \Tfrac{\Ds\delta}{\Ds\varepsilon}
  \hspace{1.5em}\forall\,\varepsilon>0                    &
\end{eqnarray}\endgroup                                        \end{Lem}
\begin{Bew}\enskip      
Choose any dominating positive measure~$\mu$, and
densities~$p$,~$q$ such that $dP=p\,d\mu$ and $ dQ=q\,d\mu $.
Then $ d \nu_c=(q-c \hspace{4\ei}p)\,d\mu$ and,
by Rudin~(1974; Theorem~6.13),  $ d|\nu_c|=|q-c \hspace{4\ei}p|\,d\mu$.
Since $ \itg {(\tau^*-\tau)\,d \nu_c} \le (1+c)\hspace{6\ei}\delta $
by~(\ref{e.x.lpt*t}), and $ (\tau^*-\tau)(q-cp)\ge0 $ a.e.$\mu$,
we conclude that
\begin{equation}
    \itg |\tau^*-\tau|\,d |\nu_c| =
    \itg {(\tau^*-\tau)\,d \nu_c}  \le (1+c)\hspace{6\ei}\delta
\end{equation}
Via the Chebyshev--Markov inequality, (\ref{e.x.nct*t}) follows.
\qed\end{Bew}
\paragraph{Acknowledgement}
I thank P.~Ruckdeschel for the numerical computations.

\vfill  
\settowidth{\breit}{\sc  e-mail:
                    \rm helmut.rieder@uni-bayreuth.de}
\hfill
\parbox{\breit}{\sc 
                    Department of Mathematics\\
                    University of Bayreuth, NW~II\\
                    D-95440 Bayreuth, Germany\\   e-mail:
                    {\rm helmut.rieder@uni-bayreuth.de}}
\end{document}